\documentclass[11pt,a4paper]{amsart}
\usepackage[T1]{fontenc}
\usepackage{lmodern}
\usepackage{amsmath,amssymb,amsthm,mathtools}
\usepackage[margin=28mm]{geometry}
\usepackage{microtype}
\usepackage{enumitem}
\usepackage[hidelinks]{hyperref}
\usepackage{bookmark}
\hypersetup{
 pdftitle={Quaternionic contact-symplectic pairs and contact pairs: symplectization and splitting},
 pdfsubject={Revised mathematical draft: definitions, reconstruction, splitting and examples},
 pdfkeywords={contact pairs, contact-symplectic pairs, hyper-Kahler geometry, 3-Sasakian geometry, symplectization}
}
\setlist[enumerate]{label=\textup{(\roman*)},leftmargin=2.6em}
\numberwithin{equation}{section}
\newtheorem{theorem}{Theorem}[section]
\newtheorem{proposition}[theorem]{Proposition}
\newtheorem{lemma}[theorem]{Lemma}

\theoremstyle{definition}
\newtheorem{definition}[theorem]{Definition}
\newtheorem{example}[theorem]{Example}
\newtheorem{remark}[theorem]{Remark}
\newcommand{\R}{\mathbb R}
\newcommand{\Z}{\mathbb Z}
\newcommand{\HH}{\mathbb H}
\newcommand{\Id}{\operatorname{Id}}
\newcommand{\rank}{\operatorname{rank}}
\newcommand{\Ric}{\operatorname{Ric}}

\newcommand{\Sp}{\operatorname{Sp}}
\newcommand{\Lie}{\mathcal L}
\newcommand{\E}{\mathcal E}
\newcommand{\F}{\mathcal F}
\newcommand{\charac}{\mathcal C}

\DeclareMathOperator{\pr}{pr}
\title{Quaternionic geometry of symplectic and contact pairs}

\author{Gianluca Bande}
\address{Dipartimento di Matematica e Informatica, Università di Cagliari, Via Ospedale 72, Cagliari, Italia}
\email{gbande@unica.it}

\author{Beniamino Cappelletti-Montano}
\address{Dipartimento di Matematica e Informatica, Università di Cagliari, Via Ospedale 72, Cagliari, Italia}
\email{b.cappellettimontano@unica.it}
\date{}

\begin{document}

\begin{abstract}
We introduce quaternionic analogues of symplectic pairs,
contact--symplectic pairs and contact pairs. A hyper symplectic
pair consists of three symplectic pairs whose sums form a
positive quaternionic triple, while the contact analogues
are defined through single or double symplectization.
We prove that radial invariance of the induced complex
structures is equivalent to the coincidence of the
characteristic distributions, yielding a local product
decomposition. Under a natural normalization, hyper
contact--symplectic pairs are locally products of a
$3$-Sasakian and a hyper-K\"ahler manifold and induce
canonical $3$-quasi-Sasakian structures, whereas hyper
contact pairs are locally products of two $3$-Sasakian
manifolds. We construct compact examples for which these
local decompositions are not global. Finally, compact
normalized hyper contact pairs have finite fundamental
group and are finite quotients of products of compact
$3$-Sasakian manifolds, in contrast with normalized hyper
contact--symplectic pairs, which may have infinite
fundamental group.
\end{abstract}


\thanks{This work was partially supported by
Fondazione di Sardegna under project ProBiki. The first author is member of the Italian National Group G.N.S.A.G.A. of INdAM. The second author acknowledge the support from GeoStat Lab}

\maketitle
\begin{center}
\itshape Dedicated to the memory of David E. Blair,\\
a mentor and a friend.
\end{center}

\section{L'introduzione}\label{sec:intro}

Symplectic pairs are pairs of closed $2$-forms of constant ranks
whose kernel foliations are complementary; each form restricts to a
symplectic form on the leaves of the kernel of the other~\cite{BK06}.
Related structures are contact--symplectic
pairs~\cite{Bande03}, consisting of a $1$-form and a closed $2$-form,
and contact pairs~\cite{BH05,BH10}, consisting of two $1$-forms. These
structures are related by symplectization: if $(\alpha,\omega)$ is a
contact--symplectic pair on $M$, then
$\bigl(d(r^2\alpha/2),\omega\bigr)$ is a symplectic pair on
$M\times\R_{>0}$, and if $(\alpha,\beta)$ is a contact pair on $M$,
then $\bigl(d(r^2\alpha/2),d(s^2\beta/2)\bigr)$ is a symplectic pair
on $M\times\R_{>0}^2$ (Lemmas~\ref{lem:csplift}
and~\ref{lem:cplift}).

In this paper we study quaternionic versions of these structures. Our
starting point is the description of hyper-K\"ahler structures by
triples of symplectic forms: by a result of Hitchin~\cite{Hitchin87},
three closed $2$-forms that form a positive quaternionic triple at
every point determine a hyper-K\"ahler structure, whose complex
structures and metric are reconstructed algebraically from the forms
(Lemmas~\ref{lem:algebra} and~\ref{lem:closed}). Accordingly, a
\emph{hyper symplectic pair} is a triple of symplectic pairs
$(\omega_i,\eta_i)$ such that the sums $\omega_i+\eta_i$ form a
positive quaternionic triple (Definition~\ref{def:sp}).
\emph{Hyper contact--symplectic pairs} and \emph{hyper contact pairs}
are defined by requiring that the lifted pairs above form hyper
symplectic pairs on $M\times\R_{>0}$ and on $M\times\R_{>0}^2$,
respectively (Definitions~\ref{def:csp} and~\ref{def:cp}). A manifold
carrying a hyper contact--symplectic pair has dimension congruent to
$3$ modulo $4$, and one carrying a hyper contact pair has dimension
congruent to $2$ modulo $4$. In both cases the tangent bundle,
stabilized by a trivial bundle of rank one or two, carries a
quaternionic Hermitian structure (Propositions~\ref{prop:cspspin}
and~\ref{prop:cpspin}).

The three symplectic pairs of a hyper symplectic pair need not have
common characteristic distributions
(Example~\ref{ex:spdifferentkernels}). Our main tool,
Theorem~\ref{thm:radial}, concerns a positive quaternionic triple
$\Omega_i=A_i+B_i$ in which the forms $A_i$ and $B_i$ are closed and
have complementary kernels, together with a vector field $R$ such
that $\Lie_RA_i=2A_i$ and $\Lie_RB_i=0$. It expresses $\Lie_RI_i$ in
terms of the projections onto $\ker B_i$ along $\ker A_i$, and shows
that the reconstructed complex structures are invariant under $R$ if
and only if the kernels do not depend on $i$. In that case the
manifold is locally the Riemannian product of two hyper-K\"ahler
manifolds. The part of Theorem~\ref{thm:radial} that does not
involve $R$ (Remark~\ref{osservazione-R}) gives
Theorem~\ref{thm:spalignment} for hyper symplectic pairs. Applied to
the symplectizations, with $R=r\partial_r$, it shows that a hyper
contact--symplectic pair or a hyper contact pair is \emph{radially
invariant}, that is, the complex structures of its symplectization
are invariant under $R$, if and only if the characteristic
distributions of the three pairs coincide
(Theorems~\ref{thm:cspalignment} and~\ref{thm:cpalignment}). In this
case $M$ splits locally as a product, and the forms are pulled back
from the factors.

Our main local results concern \emph{normalized} hyper pairs, that
is, radially invariant ones whose lifted metric satisfies
$\widehat g(R,\cdot)=r\,dr$, and also $\widehat g(S,\cdot)=s\,ds$ with
$S=s\partial_s$ in the case of contact pairs. Combining the local
splitting with a characterization of hyper-K\"ahler cones
(Lemma~\ref{lem:cone}), which relies on Kashiwada's
theorem~\cite{Kashiwada01}, we prove that a normalized hyper
contact--symplectic pair is locally, with its canonical metric, the
Riemannian product of a $3$-Sasakian manifold and a hyper-K\"ahler
manifold (Theorem~\ref{thm:cspnormalized}). In particular $M$ carries
a canonical $3$-quasi-Sasakian structure, with fundamental forms
$\frac12d\alpha_i+\omega_i$. Conversely, every $3$-quasi-Sasakian
manifold gives rise to a hyper contact--symplectic pair
(Remark~\ref{rem:3qs}), in agreement with the local description of
$3$-quasi-Sasakian manifolds in \cite{CD08} and \cite{CD09}. Similarly, a normalized
hyper contact pair is locally the Riemannian product of two
$3$-Sasakian manifolds (Theorem~\ref{thm:cpnormalized}), and
Proposition~\ref{prop:Ftensors} describes the structure tensors that
it induces on $M$.

These local product decompositions need not be global. Using
foliated bundles~\cite{CN,C}, we construct compact manifolds with
normalized hyper contact--symplectic pairs
(Proposition~\ref{prop:cspfoliated} and Examples~\ref{ex:cspfinite}
and~\ref{ex:cspinfinite}) and with normalized hyper contact pairs
(Proposition~\ref{prop:cpfoliated} and Example~\ref{ex:cpspheres})
whose canonical splitting is not a global product. The two cases
behave differently. For a normalized hyper contact pair both local
factors have positive Ricci curvature; hence, on a compact manifold,
the fundamental group is finite and the manifold is a finite quotient
of a product of two compact simply connected $3$-Sasakian manifolds
(Proposition~\ref{prop:cpcompact}). In contrast,
Example~\ref{ex:cspinfinite} gives compact manifolds with normalized
hyper contact--symplectic pairs whose fundamental group is
$\Z^{4m}$ and whose characteristic leaves are not compact
(Remark~\ref{rem:cpvscsp}).


\subsection*{Conventions}\label{subsec:conventions}
All manifolds, forms, and distributions are smooth. The exterior derivative is defined by
\[
 d\alpha(X,Y)=X(\alpha(Y))-Y(\alpha(X))-\alpha([X,Y]).
\]
It is related to the interior product and the Lie derivative by
\[ \Lie_X=d\iota_X+\iota_Xd. \]
For a $2$-form $\tau$, write $\tau^\flat X=\iota_X\tau$; when $\tau$ is nondegenerate, $\tau^\sharp=(\tau^\flat)^{-1}$. Powers of forms denote wedge powers.  
An index triple $(i,j,k)$ is cyclic when it is an even permutation of $(1,2,3)$.
For a Hermitian structure we use $\Omega_i(X,Y)=g(X,I_iY)$.
In a Sasakian manifold we have $d\alpha=2g(\cdot,\phi\cdot)$ and
\begin{equation}\label{eq:sasakiconvention}
 (\nabla_X\phi)Y=g(X,Y)\xi-\alpha(Y)X.
\end{equation}

\section{Preliminary results}\label{sec:algebra}

We recall some preliminary constructions that will be used in the following sections.

\begin{lemma}\label{lem:algebra}
Let $V$ be a finite-dimensional real vector space, let
$V^*=\operatorname{Hom}_{\mathbb R}(V,\mathbb R)$ be its dual space,
and let
$\Omega_1,\Omega_2,\Omega_3\in\Lambda^2 V^*$
be nondegenerate alternating bilinear forms.
For each $i\in\{1,2,3\}$, denote by
\[
  b_i=\Omega_i^\flat:V\longrightarrow V^*,
  \qquad
  v\longmapsto \iota_v\Omega_i=\Omega_i(v,\cdot),
\]
the isomorphism induced by $\Omega_i$. Assume that,   for all distinct $i,j\in\{1,2,3\}$,
\begin{equation}\label{eq:algebraiccondition}
  b_i^{-1}\circ b_j=-b_j^{-1}\circ b_i.
\end{equation}
For each cyclic permutation $(i,j,k)$ of $(1,2,3)$, define
\begin{equation}\label{eq:reconstruction}
  I_i:=-b_j^{-1}\circ b_k:V\longrightarrow V,
  \qquad
  G_i:=b_i\circ I_i:V\longrightarrow V^*.
\end{equation}
Then
  $G_1=G_2=G_3$.
Denote their common value by $G:V\longrightarrow V^*$ and define
\begin{equation}\label{reconstructed}
  g:V\times V\longrightarrow\mathbb R,
  \qquad
  g(v,w):=\bigl(G(v)\bigr)(w).
\end{equation}
The bilinear form $g$ is symmetric and nondegenerate. Moreover,
the endomorphisms $I_1,I_2,I_3$ satisfy
\[
  I_i^2=-\operatorname{Id}_V,
  \qquad
  I_i\circ I_j=I_k=-I_j\circ I_i
\]
for every cyclic permutation $(i,j,k)$, and
  $g(I_i v,I_i w)=g(v,w)$
  for all $v,w\in V$
   and $i\in\{1,2,3\}$.
Finally, for all $v,w\in V$ and $i\in\{1,2,3\}$,
\begin{equation}\label{eq:metricreconstruction}
  \Omega_i(v,w)=g(v,I_i w),
  \qquad
  g(v,w)=-\Omega_i(v,I_i w).
\end{equation}
\end{lemma}

\begin{proof}
Fix a cyclic permutation $(i,j,k)$ of $(1,2,3)$.
Using \eqref{eq:algebraiccondition}, we obtain
\[
  I_i^2
  =(b_j^{-1}\circ b_k)\circ(b_j^{-1}\circ b_k)
  =-(b_j^{-1}\circ b_k)\circ(b_k^{-1}\circ b_j)
  =-\operatorname{Id}_V.
\]
Furthermore,
\[
  I_i\circ I_j
  =b_j^{-1}\circ b_k\circ b_k^{-1}\circ b_i
  =b_j^{-1}\circ b_i
  =-b_i^{-1}\circ b_j
  =I_k.
\]
Applying this identity cyclically also gives
  $I_j\circ I_i
  =I_j\circ(I_j\circ I_k)
  =-I_k$.
Thus $I_1,I_2,I_3$ satisfy the quaternionic identities.
Composing \eqref{eq:algebraiccondition} on the left by $b_i$
and on the right by $b_i^{-1}$ gives
  $b_i\circ b_j^{-1}=-b_j\circ b_i^{-1}$.
Hence,
\begin{align*}
  G_i=-b_i\circ b_j^{-1}\circ b_k=b_j\circ b_i^{-1}\circ b_k=-b_j\circ b_k^{-1}\circ b_i=G_j.
\end{align*}
Set $G:=G_1=G_2=G_3$.
Since $G=b_i\circ I_i$ is a composition of isomorphisms,
the associated bilinear form
  $g(v,w)=\bigl(G(v)\bigr)(w)=\Omega_i(I_i v,w)$
is nondegenerate.
We next prove that $g$ is symmetric.
The definitions and \eqref{eq:algebraiccondition} imply
\[
  b_i\circ I_j=b_k,
  \qquad
  b_i\circ I_k=-b_j.
\]
Since $\Omega_j$ and $\Omega_k$ are alternating, these identities
yield
  $\Omega_i(I_a v,w)=\Omega_i(v,I_a w)$
  for all $a\in\{j,k\}$
   $v,w\in V$.
Therefore
\[
  g(v,w)
  =-\Omega_i(v,I_i w)
  =\Omega_i(I_i w,v)
  =g(w,v),
\]
which proves symmetry and the second identity
in \eqref{eq:metricreconstruction}.
Replacing $w$ by $I_i w$ in that identity gives
\[
  g(v,I_i w)
  =-\Omega_i(v,I_i^2w)
  =\Omega_i(v,w),
\]
which is the first identity in
\eqref{eq:metricreconstruction}.
Finally,
\[
  g(I_i v,I_i w)
  =\Omega_i(I_i^2v,I_i w)
  =-\Omega_i(v,I_i w)
  =g(v,w).
\]
\end{proof}

The metric $g$ defined by \eqref{reconstructed} will be called the \emph{reconstructed metric} from the nondegenerate alternating bilinear forms $\Omega_1, \Omega_2, \Omega_3$ satisfying \eqref{eq:algebraiccondition}.

\begin{definition}\label{def:positive}
A triple $(\Omega_1,\Omega_2,\Omega_3)$ is a \emph{positive quaternionic triple} if it satisfies Lemma~\ref{lem:algebra} and the reconstructed metric $g$ is positive definite.
\end{definition}

For forms on a manifold the above conditions are imposed at every point.
The following result is a reformulation in our language of a well-known result of Hitchin (\cite{Hitchin87}).

\begin{lemma}\label{lem:closed}
Let $M$ be a smooth real manifold, and let
 $\Omega_1,\Omega_2,\Omega_3\in\Omega^2(M)$
be smooth closed $2$-forms such that, at every point $p\in M$,
the triple
  $\left((\Omega_1)_p,(\Omega_2)_p,(\Omega_3)_p\right)$
is a positive quaternionic triple on $T_pM$.
Let $I_1,I_2,I_3:TM\longrightarrow TM$ and $g$ be obtained
by applying Lemma~\ref{lem:algebra} pointwise, that is
\[
  I_i=-(\Omega_j^\flat)^{-1}\circ\Omega_k^\flat
\]
for every cyclic permutation $(i,j,k)$ of $(1,2,3)$, and
\[
  g_p(v,w)=-(\Omega_i)_p\bigl(v,(I_i)_p w\bigr),
  \qquad v,w\in T_pM,
\]
independently of $i\in\left\{1,2,3\right\}$.
Then $(g,I_1,I_2,I_3)$ is a hyper-K\"ahler structure on $M$.
\end{lemma}

\begin{remark}
In the literature, the endomorphisms $I_j$ are called \textit{Recursion Operators}. In \cite{BK01} it was shown that, if a Recursion Operator squares to minus the identity, actually it is authomatically integrable.     
\end{remark}

\begin{theorem}\label{thm:radial}
Let $\Omega_i=A_i+B_i$ be a positive quaternionic triple on a manifold $M$. Assume that each $A_i,B_i$ is closed, has constant rank, and satisfies
\begin{equation}\label{eq:complementaryAB}
 \ker A_i\oplus\ker B_i=T M.
\end{equation}
Let $R$ be a vector field such that
\begin{equation}\label{eq:radialweights}
 \Lie_RA_i=2A_i,\qquad \Lie_RB_i=0.
\end{equation}
Denote by $I_i, g$ the reconstructed structure according to Lemma \ref{lem:closed} and set
 $\Pi_i=(\Omega_i^\flat)^{-1}A_i^\flat$.
Then $\Pi_i$ is the projection onto $\ker B_i$ along $\ker A_i$, and
\begin{equation}\label{eq:radialobstruction}
 \Lie_RI_i=2(I_i\Pi_k-\Pi_jI_i)
\end{equation}
for any cyclic permutation $(i,j,k)$.
Moreover, the following conditions are equivalent:
\begin{enumerate}
\item $\Lie_RI_i=0$ for $i=1,2,3$;
\item $\Pi_1=\Pi_2=\Pi_3$;
\item both $\ker A_i$ and $\ker B_i$ are independent of $i$.
\end{enumerate}
If these conditions hold, the common projector $\Pi$ is $ g$-orthogonal, commutes with all $I_i$, and is parallel. Its image and kernel have dimensions divisible by $4$. Locally, $M$ is the Riemannian product of two hyper-K\"ahler manifolds, and $A_i,B_i$ are the respective fundamental forms.
\end{theorem}
\begin{proof}
We first prove the assertion concerning the projectors.
Fix $i\in\{1,2,3\}$ and $p\in M$. We work on the
vector space $V=T_p  M$. For brevity, write
  $A=(A_i)_p$,  $B=(B_i)_p$,
  $\Omega=(\Omega_i)_p=A+B$,
and
\[
  P=(\Pi_i)_p
   =(\Omega^\flat)^{-1}\circ A^\flat:V\longrightarrow V.
\]
By \eqref{eq:complementaryAB}, every $v\in V$ admits
a unique decomposition
  $v=v_A+v_B$,
with
  $v_A\in\ker A$ and $v_B\in\ker B$.
For every $w\in V$, the definition of $P$ implies
\begin{align*}
  \Omega(Pv,w)
  =A(v,w)=A(v_B,w)=A(v_B,w)+B(v_B,w)=\Omega(v_B,w),
\end{align*}
where we have used  $v_A\in\ker A$ and $v_B\in\ker B$.
Hence,
  $\Omega(Pv-v_B,w)=0$
  for every $w\in V$.
Since $\Omega$ is nondegenerate, it follows that
  $P(v_A+v_B)=v_B$.
Thus $P$ is precisely the projection onto $\ker B$
along $\ker A$. In particular,
  $P|_{\ker A}=0$ and
  $P|_{\ker B}=\operatorname{Id}_{\ker B}$,
and therefore
\[
  P^2=P,\qquad
  \ker P=\ker A,\qquad
  \operatorname{im}P=\ker B.
\]
As $p$ was arbitrary, this proves that $\Pi_i$ is the
smooth bundle projection onto $\ker B_i$ along $\ker A_i$.
In particular,
\[
  \Pi_i^2=\Pi_i,\qquad
  \ker\Pi_i=\ker A_i,\qquad
  \operatorname{im}\Pi_i=\ker B_i.
\]
These identities also prove the equivalence of
\textup{(ii)} and \textup{(iii)}.
Next, differentiating $I_i=-(\Omega_j^\flat)^{-1}\Omega_k^\flat$ and using $\Lie_R\Omega_a=2A_a$ give
\begin{align*}
 \Lie_RI_i=2(\Omega_j^\flat)^{-1}A_j^\flat
          (\Omega_j^\flat)^{-1}\Omega_k^\flat
   -2(\Omega_j^\flat)^{-1}A_k^\flat=-2\Pi_jI_i+2I_i\Pi_k,
\end{align*}
which proves \eqref{eq:radialobstruction}.
Now, assume (i) and set $K=\frac12\Lie_R  g$, a symmetric bilinear form. Differentiating $  g(I_aX,I_aY)=  g(X,Y)$ shows that $K$ is invariant under every $I_a$. Thus the endomorphism $H$ defined by
 $K(X,Y)= g(HX,Y)$
commutes with every $I_a$. Differentiating $\Omega_i(X,Y)=  g(X,I_iY)$ now gives
$A_i(X,Y)=K(X,I_iY)=\Omega_i(HX,Y)$.
Consequently $\Pi_i=H$ for every $i$, proving (ii).
Conversely, assume \textup{(iii)} and set
  $\E:=\ker B_i$, $\F:=\ker A_i$,
which are independent of $i$ by hypothesis. Thus
  $T M=\E\oplus\F$.
Denote by $\Pi$ the projection onto $\E$ along $\F$.
By the first part of the proof,
  $\Pi_1=\Pi_2=\Pi_3=\Pi$.
We first examine this decomposition pointwise.
Fix $p\in M$ and write
\[
  V:=T_p M,\qquad E:=\E_p,\qquad F:=\F_p,
  \qquad \omega_i:=(\Omega_i)_p.
\]
For every $v\in E$ and $w\in F$, we have
  $\omega_i(v,w)
  =(A_i)_p(v,w)+(B_i)_p(v,w)=0$,
because $w\in\ker(A_i)_p$ and $v\in\ker(B_i)_p$.
Hence, if
  $\omega_i^E:=\omega_i|_{E\times E}$,
and 
  $\omega_i^F:=\omega_i|_{F\times F}$,
then
  $\omega_i(v_E+v_F,w_E+w_F)
=\omega_i^E(v_E,w_E)+\omega_i^F(v_F,w_F)$
for all $v_E,w_E\in E$ and $v_F,w_F\in F$.
Both restricted forms are nondegenerate.
Indeed, if $v\in E$ satisfies
  $\omega_i^E(v,w)=0$ for every $w\in E$,
then, by the vanishing of the mixed terms, $\omega_i(v,w)=0$ also for every $w\in F$.
Thus $\omega_i(v,\cdot)=0$ on $V$, and the nondegeneracy
of $\omega_i$ implies $v=0$.
The same argument applies to $\omega_i^F$.
The splitting $V=E\oplus F$ induces the identification
\[
  V^*\longrightarrow E^*\oplus F^*,
  \qquad
  \lambda\longmapsto(\lambda|_E,\lambda|_F).
\]
This identification uses only the direct sum decomposition,
not the metric. With respect to these decomposition, set
\[
  b_i^E:=(\omega_i^E)^\flat:E\longrightarrow E^*,
  \qquad
  b_i^F:=(\omega_i^F)^\flat:F\longrightarrow F^*.
\]
The preceding identities give
\[
  \omega_i^\flat=
  \begin{pmatrix}
    b_i^E & 0\\
    0 & b_i^F
  \end{pmatrix},
  \qquad
  (\omega_i^\flat)^{-1}=
  \begin{pmatrix}
    (b_i^E)^{-1} & 0\\
    0 & (b_i^F)^{-1}
  \end{pmatrix}.
\]
It follows that, for each cyclic permutation $(i,j,k)$,
\begin{equation}\label{blocchi}
  (I_i)_p
  =-(\omega_j^\flat)^{-1}\circ\omega_k^\flat
  =
  \begin{pmatrix}
    -(b_j^E)^{-1}\circ b_k^E & 0\\
    0 & -(b_j^F)^{-1}\circ b_k^F
  \end{pmatrix}.
\end{equation}
In particular,
  $(I_i)_p(E)\subseteq E$,
  and
  $(I_i)_p(F)\subseteq F$.
Since $p$ was arbitrary, each $I_i$ preserves both
distributions $\E$ and $\F$.
Moreover, $(\Pi)_p$ acts as the identity on $E$ and
vanishes on $F$. Thus, for $v=v_E+v_F$,
\[
  (I_i)_p(\Pi)_p v
  =(I_i)_p v_E
  =(\Pi)_p(I_i)_p v.
\]
Therefore
  $I_i\Pi=\Pi I_i$.
Substituting $\Pi_j=\Pi_k=\Pi$ into
\eqref{eq:radialobstruction}, we obtain
  $\mathcal L_R I_i
  =2(I_i\Pi-\Pi I_i)=0$.
This proves \textup{(i)}.
We next verify the metric assertions.
For $v\in E$ and $w\in F$, we have
  $g_p(v,w)
  =-\omega_i\bigl(v,(I_i)_p w\bigr)=0$,
since $(I_i)_p w\in F$ and $\omega_i(E,F)=0$.
Hence $\E$ and $\F$ are $g$-orthogonal.
As $\Pi$ is the identity on $\E$ and vanishes on $\F$,
it is therefore the $g$-orthogonal projection onto $\E$.
More explicitly, write
\[
  I_i^E:=(I_i)_p|_E,\qquad I_i^F:=(I_i)_p|_F.
\]
The block formulas \eqref{blocchi} show that these are precisely
the endomorphisms reconstructed from the restricted triples
$(\omega_1^E,\omega_2^E,\omega_3^E)$ and
$(\omega_1^F,\omega_2^F,\omega_3^F)$, respectively.
The corresponding reconstructed bilinear forms satisfy
\[
  -\omega_i^E(v,I_i^Ew)=g_p(v,w)
  \qquad\text{for }v,w\in E,
\]
and
\[
  -\omega_i^F(v,I_i^Fw)=g_p(v,w)
  \qquad\text{for }v,w\in F.
\]
Thus the metrics reconstructed on the two blocks are
exactly $g_p|_{E\times E}$ and $g_p|_{F\times F}$.
In particular, they are positive definite.
Finally, the quaternionic identities for the $I_i$
restrict to both $E$ and $F$, because these subspaces
are invariant under every $(I_i)_p$.
Consequently, each block is a quaternionic Hermitian
vector space. Therefore
  $\dim E\equiv\dim F\equiv0\pmod4$.
It remains to prove the local product assertion and the parallelism
of $\Pi$. Since the forms are closed and have constant rank, their
kernel distributions $\E$ and $\F$ are integrable. Their
complementarity gives a local product of leaves. For
$X\in\Gamma(\E)$ and $Y\in\Gamma(\F)$, Cartan's formula gives
\[
  \Lie_YA_i=d(\iota_YA_i)+\iota_YdA_i=0,
  \qquad
  \Lie_XB_i=d(\iota_XB_i)+\iota_XdB_i=0.
\]
Together with $\iota_YA_i=0$ and $\iota_XB_i=0$, these identities
show that $A_i$ is pulled back from the $\E$-factor and $B_i$ from
the $\F$-factor. The blockwise
reconstruction established above therefore gives a metric on each
factor, and $g$ is their Riemannian product. The restricted triples
are positive and closed, so the two factors are hyper-K\"ahler by
Lemma~\ref{lem:closed}. The Levi-Civita connection of a product
metric preserves both tangent summands. Since $\Pi$ acts as the
identity on one summand and as zero on the other, $\nabla\Pi=0$.
\end{proof}

\begin{remark}\label{osservazione-R}
Note that the statements in Theorem \ref{thm:radial}
concerning the projectors, the equivalence
of \textup{(ii)} and \textup{(iii)}, and the geometric
conclusions under either of these two conditions remain
valid without assuming the existence of $R$.
The additional hypothesis \eqref{eq:radialweights}
is needed only for \eqref{eq:radialobstruction}
and for the equivalence with \textup{(i)}.
\end{remark}

\section{Hyper symplectic pairs}\label{sec:sp}

Recall that a \emph{symplectic pair} $(\omega,\eta)$ of type $(p,q)$
on a manifold of dimension $2p+2q$ consists of two closed $2$-forms such
that
\begin{equation}\label{eq:spordinary}
  \omega^p\wedge \eta^q\neq0,
  \qquad \omega^{p+1}=0,
  \qquad \eta^{q+1}=0.
\end{equation}
Then $\omega$ and $\eta$ have constant ranks $2p$ and $2q$, respectively, and
\[
  TM=\ker \eta\oplus\ker \omega.
\]
Each form restricts to a symplectic form on the leaves of the kernel foliation
of the other.
For every symplectic pair, the form $\Omega=\omega+\eta$ is symplectic.
Indeed, it is closed and
\begin{equation}\label{somma-simplettica}
  (\omega+\eta)^{p+q}=\binom{p+q}{p}\omega^p\wedge \eta^q\neq0.
\end{equation}
We refer to \cite{BK06} for more details on symplectic pairs.

\begin{definition}\label{def:sp}
A \emph{hyper symplectic pair} on $M$ is a triple of symplectic pairs $(\omega_1,\eta_1)$, $(\omega_2,\eta_2)$, $(\omega_3,\eta_3)$
such that the $2$-forms
  $\Omega_i:=\omega_i+\eta_i$, $i\in\{1,2,3\}$,
form a positive quaternionic triple on $M$.
The complex structures and metric reconstructed from the $\Omega_i$, according to  Lemma \ref{lem:algebra}, are
\begin{equation}\label{eq:spreconstruction}
  I_i=-(\Omega_j^\flat)^{-1}\circ\Omega_k^\flat,
  \qquad
  g(X,Y)=-\Omega_i(X,I_iY),
\end{equation}
for cyclic $(i,j,k)$. We call $g$ the \emph{canonical metric} of the
hyper symplectic pair.
\end{definition}

\begin{proposition}\label{prop:sphyperkahler}
A hyper symplectic pair determines a hyper-K\"ahler structure
$(g,I_1,I_2,I_3)$ on $M$. In particular,
  $\dim M\equiv0\pmod4$.
\end{proposition}
\begin{proof}
 This is an immediate application of Lemmas~\ref{lem:algebra}
and~\ref{lem:closed}, since $d\Omega_i=d\omega_i+d\eta_i=0$.   
\end{proof}

\begin{remark}\label{rem:spdata}
The definition specifies more than a hyper-K\"ahler structure:
it includes three decompositions of its fundamental forms into
symplectic pairs. In general, the three pairs need not have common
characteristic distributions. 
\end{remark}

For each $i\in\left\{1,2,3\right\}$, we set
  $\E_i:=\ker\eta_i$, $\F_i:=\ker\omega_i$,
   $\Pi_i:=(\Omega_i^\flat)^{-1}\circ\omega_i^\flat$.
The pointwise computation at the beginning of the proof of
Theorem~\ref{thm:radial} shows that $\Pi_i$ is the projection onto
$\E_i$ along $\F_i$ (note that the computation uses only the complementarity
of the two kernels).

\begin{theorem}\label{thm:spalignment}
Let $(\omega_i,\eta_i)_{i=1,2,3}$ be a hyper symplectic pair.
The following conditions are equivalent:
\begin{enumerate}
\item
The three symplectic pairs have common characteristic distributions, namely
\begin{equation}\label{eq:spalignment}
  \ker \eta_1=\ker \eta_2=\ker \eta_3=\E,
  \qquad
  \ker \omega_1=\ker \omega_2=\ker \omega_3=\F.
\end{equation}
\item
The projectors satisfy $\Pi_1=\Pi_2=\Pi_3$.
\end{enumerate}
If these conditions hold, all three pairs have the same type
$(2l,2m)$, for some $l,m\geq0$, and
\[
  TM=\E\oplus\F,
  \qquad \rank\E=4l,
  \qquad \rank\F=4m.
\]
The common projector $\Pi$ is $g$-orthogonal, commutes with every
$I_i$, and is parallel with respect to the Levi-Civita connection of $g$.
Locally there is a Riemannian product
\begin{equation}\label{eq:spproduct}
  (M,g)\cong(P^{4l},g_P)\times(Q^{4m},g_Q)
\end{equation}
of hyper-K\"ahler manifolds such that
  $\omega_i=\pr_P^*\Omega_i^P$
and
  $\eta_i=\pr_Q^*\Omega_i^Q$,
where $\Omega_i^P$ and $\Omega_i^Q$ are the respective fundamental
forms. Conversely, globally defined forms admitting the local description above determine a hyper symplectic pair with common characteristic distributions.
\end{theorem}

\begin{proof}
Let $(p_i,q_i)$ be the type of the symplectic pair
$(\omega_i,\eta_i)$, and write
  $\Omega_i=\omega_i+\eta_i$,
  $\Pi_i=(\Omega_i^\flat)^{-1}\circ\omega_i^\flat$.
By definition, the forms $\omega_i,\eta_i$ are closed,
have constant rank, and satisfy
\[
  TM=\ker\eta_i\oplus\ker\omega_i.
\]
Moreover, $(\Omega_1,\Omega_2,\Omega_3)$ is a positive
quaternionic triple.
We apply the part of Theorem~\ref{thm:radial} that
does not require the existence of $R$ (see Remark \ref{osservazione-R}), with
  $A_i=\omega_i$ and $B_i=\eta_i$.
It follows that $\Pi_i$ is the projection onto
$\ker\eta_i$ along $\ker\omega_i$. In particular,
\begin{equation}\label{projections}
  \Pi_i^2=\Pi_i,\qquad
  \operatorname{im}\Pi_i=\ker\eta_i,\qquad
  \ker\Pi_i=\ker\omega_i.
\end{equation}
Thus, the equality of the three projectors is equivalent
to the equality of their respective images and kernels,
which proves the equivalence of \textup{(i)}
and \textup{(ii)}.
Assume that these equivalent conditions hold, and put
\[
  \E=\ker\eta_i,\qquad
  \F=\ker\omega_i,\qquad
  \Pi=\Pi_1=\Pi_2=\Pi_3.
\]
By Theorem~\ref{thm:radial} we have 
\[
  \E\perp_g\F,\qquad
  I_i\Pi=\Pi I_i,\qquad
  \nabla\Pi=0,
\]
and the ranks of $\E$ and $\F$ are divisible
by $4$. Write
  $\operatorname{rank}\E=4l$ and
  $\operatorname{rank}\F=4m$
for some integers $l,m\geq0$.
Since the $i$th symplectic pair has type $(p_i,q_i)$, we have
 $ \operatorname{rank}\omega_i=2p_i$ and
  $\operatorname{rank}\eta_i=2q_i$.
Using their kernels and the complementary splitting
$TM=\E\oplus\F$, we have 
\[
  2p_i=\dim M-\operatorname{rank}\F=4l,
  \qquad
  2q_i=\dim M-\operatorname{rank}\E=4m.
\]
Hence all three pairs have type $(2l,2m)$.
Finally, by
Theorem~\ref{thm:radial}  we have, locally, 
\[
(M,g)\cong(P^{4l},g_P)\times(Q^{4m},g_Q)
\]
with both factors hyper-K\"ahler and
  $\omega_i=\pr_P^*\Omega_i^P$, 
  $\eta_i=\pr_Q^*\Omega_i^Q$.
Conversely, suppose the globally defined forms admit
such a local product description.
The forms $\omega_i$ and $\eta_i$ are closed and satisf
  $\omega_i^{2l+1}=0$,
  $\eta_i^{2m+1}=0$,
whereas
\[
  \omega_i^{2l}\wedge\eta_i^{2m}
  =
  \pr_P^*\bigl((\Omega_i^P)^{2l}\bigr)
  \wedge
  \pr_Q^*\bigl((\Omega_i^Q)^{2m}\bigr)
\]
is nowhere vanishing.
Thus each $(\omega_i,\eta_i)$ is a symplectic pair
of type $(2l,2m)$.
If $J_i$ and $K_i$ are the complex structures
on $P$ and $Q$, respectively, then
  $I_i=J_i\oplus K_i$
satisfy the quaternionic identities and are orthogonal
for the product metric $g_P+g_Q$.
Their fundamental forms are
  $\Omega_i=\omega_i+\eta_i$.
Hence these sums form a positive quaternionic triple,
so the given forms define a hyper symplectic pair.
Moreover, on every such neighborhood,
  $\ker\eta_i=TP$,
  $\ker\omega_i=TQ$
for all $i$, proving that the characteristic
distributions are common.
Since the forms are globally defined, these local
verifications establish the asserted global properties.
On the intersection of two such product neighborhoods,
the forms $\Omega_i=\omega_i+\eta_i$ are the restrictions
of the same globally defined two-forms. Since the reconstructed
endomorphisms are given pointwise by
  $I_i=-(\Omega_j^\flat)^{-1}\circ\Omega_k^\flat$,
they coincide on the overlap. The corresponding metrics
also coincide, because
  $g(X,Y)=-\Omega_i(X,I_iY)$.
Hence the locally reconstructed tensors glue to the globally
defined structure determined by the forms $\Omega_i$.
\end{proof}

\begin{example}\label{ex:spproduct}
Let $P^{4l}$ and $Q^{4m}$ be hyper-K\"ahler manifolds, with fundamental
forms $\Omega_i^P$ and $\Omega_i^Q$. On $P\times Q$, the pairs
  $(\omega_i,\eta_i)=(\pr_P^*\Omega_i^P,\pr_Q^*\Omega_i^Q)$
define a hyper symplectic pair of type $(2l,2m)$ with common
characteristic distributions and canonical metric $g_P+g_Q$.
For example, the product of flat quaternionic tori
$T^{4l}\times T^{4m}$ gives compact examples.
\end{example}

\begin{example}\label{ex:spdifferentkernels}
Common characteristic distributions are not automatic. Indeed
let $M=T^4=\R^4/\Z^4$, and let $e^0,e^1,e^2,e^3$ be the global
closed $1$-forms induced by the standard coordinate forms on
$\R^4$. Denote their dual vector fields by $e_0,e_1,e_2,e_3$.
For cyclic $(i,j,k)$ in $\{1,2,3\}$, set
\begin{equation}\label{eq:sptorusforms}
  \omega_i=e^0\wedge e^i,
  \qquad
  \eta_i=-e^j\wedge e^k.
\end{equation}
Each $(\omega_i,\eta_i)$ is a symplectic pair of type $(1,1)$, since
 $\omega_i\wedge \omega_i=\eta_i\wedge \eta_i=0$ and
  $\omega_i\wedge \eta_i=-e^0\wedge e^1\wedge e^2\wedge e^3\neq0$.
Their sums are the fundamental forms of the flat hyper-K\"ahler
structure with metric
  $g=\sum_{a=0}^3e^a\otimes e^a$
and complex structures given by
\[
  I_ie_0=-e_i,\qquad I_ie_i=e_0,
  \qquad I_ie_j=e_k,\qquad I_ie_k=-e_j.
\]
These formulas satisfy the quaternionic identities and
$\Omega_i=g(\cdot,I_i\cdot)$.
Hence \eqref{eq:sptorusforms} defines a hyper symplectic pair.
However,
$\ker\omega_i=\operatorname{span}\{e_j,e_k\}$
  and
$\ker\eta_i=\operatorname{span}\{e_0,e_i\}$,
so both kernels depend on $i\in\left\{1,2,3\right\}$.
\end{example}

\section{Hyper contact-symplectic pairs}\label{sec:csp}
A contact--symplectic pair $(\alpha,\omega)$ of type $(p,q)$ consists of a $1$-form and a closed $2$-form on a manifold of dimension $2p+2q+1$ such that
\begin{equation}\label{eq:cspordinary}
 \alpha\wedge(d\alpha)^p\wedge\omega^q\ne0,
 \qquad (d\alpha)^{p+1}=0,\qquad\omega^{q+1}=0.
\end{equation}
The characteristic distributions are the complementary distributions defined by
$\charac(\alpha)=\ker\alpha\cap\ker d\alpha$ and
  $\ker\omega$.
The Reeb vector field $Z$ is characterized by $\alpha(Z)=1$ and $\iota_Zd\alpha=\iota_Z\omega=0$.

\begin{lemma}\label{lem:csplift}
For a contact--symplectic pair of type $(p,q)$, the forms
 $A=d\!\left(\frac{r^2}{2}\alpha\right)$ and $B=\omega$
on $\widehat M=M\times\R_{>0}$ have complementary kernels
\begin{equation}\label{eq:cspkernels}
 \ker A=\charac(\alpha),\qquad
 \ker B=\ker\omega\oplus\R\partial_r.
\end{equation}
They have ranks $2p+2$ and $2q$, respectively. Thus $(A,B)$ is a
symplectic pair of type $(p+1,q)$, and $A+B$ is symplectic.
\end{lemma}
\begin{proof}
Let
  $\pi:\widehat M=M\times\mathbb R_{>0}\longrightarrow M$
be the canonical projection. We use the same notation
for the forms $\alpha$, $d\alpha$, and $\omega$ on $M$ and for
their pullbacks to $\widehat M$.
Since
  $A=d\!\left(\frac{r^2}{2}\alpha\right)$,
we have
\begin{equation}\label{eq:csp-A-explicit}
  A=r\,dr\wedge\alpha+\frac{r^2}{2}\,d\alpha.
\end{equation}
Fix a point $(x,r)\in\widehat M$. Using the canonical
identification
  $T_{(x,r)}\widehat M
  =T_xM\oplus\mathbb R\partial_r$,
write an arbitrary tangent vector as
  $v+a\partial_r$, with
  $ v\in T_xM$, $a\in\mathbb R$.
Contracting \eqref{eq:csp-A-explicit}, we obtain
\begin{align}
  \iota_{v+a\partial_r}A=
  r\,\iota_{v+a\partial_r}(dr\wedge\alpha)
  +\frac{r^2}{2}\iota_vd\alpha=
  ra\,\alpha-r\alpha(v)\,dr
  +\frac{r^2}{2}\iota_vd\alpha.
  \label{eq:csp-contraction-A}
\end{align}
We first determine $\ker A$. Suppose that
  $\iota_{v+a\partial_r}A=0$.
Evaluating \eqref{eq:csp-contraction-A} on $\partial_r$
$0=\bigl(\iota_{v+a\partial_r}A\bigr)(\partial_r)
  =-r\alpha(v)$.
Since $r>0$, it follows that
  $\alpha(v)=0$.
Let $Z$ denote the Reeb vector field of the
contact--symplectic pair $(\alpha,\omega)$, characterized by
  $\alpha(Z)=1$,
  $\iota_Zd\alpha=0$,
  $\iota_Z\omega=0$.
We regard $Z_x$ as a tangent vector to $\widehat M$ with zero
$\partial_r$-component. Evaluating
\eqref{eq:csp-contraction-A} on $Z_x$ gives
\begin{align*}
  0
=\bigl(\iota_{v+a\partial_r}A\bigr)(Z_x)
  =
  ra\,\alpha(Z_x)
  +\frac{r^2}{2}d\alpha(v,Z_x)=ra,
\end{align*}
because $\alpha(Z)=1$ and $\iota_Zd\alpha=0$.
Hence
 $a=0$.
Equation \eqref{eq:csp-contraction-A} now reduces to
  $\frac{r^2}{2}\iota_vd\alpha=0$,
and therefore
  $\iota_vd\alpha=0$.
We have proved that
  $v\in\ker\alpha\cap\ker d\alpha
  =\charac(\alpha)$.
Conversely, if
$v\in\charac(\alpha)$,
then $\alpha(v)=0$ and $\iota_vd\alpha=0$, so
\eqref{eq:csp-contraction-A} immediately gives
 $ \iota_vA=0$.
Thus
\begin{equation}\label{eq:csp-kernel-A-detail}
  \ker A=\charac(\alpha),
\end{equation}
where $\charac(\alpha)$ is regarded as a subbundle of
$T\widehat M$ contained in the $TM$-factor.
We next determine $\ker B$. Since $B=\omega$ is pulled back
from $M$, it contains no $dr$-component. Consequently,
$\iota_{v+a\partial_r}B=\iota_v\omega$.
Therefore
  $\iota_{v+a\partial_r}B=0$ if and only if
  $v\in\ker\omega$,
with no restriction on $a$. Hence
\begin{equation}\label{eq:csp-kernel-B-detail}
  \ker B
  =\ker\omega\oplus\mathbb R\partial_r.
\end{equation}
We now verify the ranks. Since $(\alpha,\omega)$ has type
$(p,q)$,
  $(d\alpha)^{p+1}=0$ and
  $\omega^{q+1}=0$,
whereas
$\alpha\wedge(d\alpha)^p\wedge\omega^q$
is nowhere vanishing. In particular,
  $(d\alpha)^p\neq0$
and
  $\omega^q\neq0$.
Thus
  $\operatorname{rank}(d\alpha)=2p$
  and
$\operatorname{rank}\omega=2q$.
For $A$, using \eqref{eq:csp-A-explicit} and 
  $(dr\wedge\alpha)^2=0$,
  $(d\alpha)^{p+1}=0$,
we get
\begin{align}
  A^{p+1}=
  \left(
 r\,dr\wedge\alpha+\frac{r^2}{2}d\alpha
  \right)^{p+1}=
  (p+1)\,
  r\,dr\wedge\alpha\wedge
  \left(\frac{r^2}{2}d\alpha\right)^p=
  \frac{(p+1)r^{2p+1}}{2^p}\,
dr\wedge\alpha\wedge(d\alpha)^p.
  \label{eq:csp-A-power}
\end{align}
Moreover,
  $A^{p+2}=0$.
The form in \eqref{eq:csp-A-power} is nonzero, since
$dr\wedge\alpha\wedge(d\alpha)^p\wedge\omega^q$
is a volume form on $\widehat M$. Hence
  $\operatorname{rank}A=2p+2$.
Since $B=\omega$, we have
  $\operatorname{rank}B=2q$.
Let us also verify explicitly that the two kernels are
complementary. If
  $u\in\ker A\cap\ker B$,
then by \eqref{eq:csp-kernel-A-detail}--%
\eqref{eq:csp-kernel-B-detail} the vector $u$ has no
$\partial_r$-component and satisfies
\[
  \alpha(u)=0,\qquad
  \iota_ud\alpha=0,\qquad
  \iota_u\omega=0.
\]
Hence,
  $\iota_u
  \bigl(
    \alpha\wedge(d\alpha)^p\wedge\omega^q
  \bigr)=0$.
Since
$\alpha\wedge(d\alpha)^p\wedge\omega^q$
is a volume form on $M$, contraction with it is injective
on tangent vectors. Hence $u=0$. Thus
  $\ker A\cap\ker B=\{0\}$.
Furthermore,
\begin{align*}
  \dim(\ker A)
  =
  \dim\widehat M-\operatorname{rank}A
   =2q, \quad
  \dim(\ker B)
  =
  \dim\widehat M-\operatorname{rank}B
   =2p+2.
\end{align*}
Since
  $2q+(2p+2)=\dim\widehat M$,
we conclude that
  $T\widehat M=\ker A\oplus\ker B$.
It remains to prove that $A+B$ is symplectic.
It is closed, since both forms are closed.
To prove nondegeneracy, set
  $n=p+q+1$,
so that $\dim\widehat M=2n$. Since
  $A^{p+2}=0$, $B^{q+1}=0$,
the only nonzero term in the binomial expansion of
$(A+B)^n$ is the one containing $A^{p+1}$ and $B^q$.
Therefore
\begin{align*}
  (A+B)^n
  &=
  \binom{n}{p+1}A^{p+1}\wedge B^q=
  \binom{p+q+1}{p+1}
  \frac{(p+1)r^{2p+1}}{2^p}\,
  dr\wedge\alpha\wedge(d\alpha)^p\wedge\omega^q.
\end{align*}
The right-hand side is nowhere vanishing. Hence
$A+B$ is nondegenerate and therefore symplectic.
\end{proof}

\begin{definition}\label{def:csp}
Let $(\alpha_i,\omega_i)$, $i=1,2,3$, be three contact--symplectic
pairs on $M$, and let
  $\pi:\widehat M=M\times\R_{>0}\longrightarrow M$
be the projection. The triple is called a
\emph{hyper contact--symplectic pair} if the pairs
\[
  \widehat\omega_i
  :=d\!\left(\frac{r^2}{2}\pi^*\alpha_i\right),
  \qquad
  \widehat\eta_i:=\pi^*\omega_i,
  \qquad i=1,2,3,
\]
form a hyper symplectic pair on $\widehat M$ in the sense of
Definition~\ref{def:sp}.
\end{definition}

Note that by Lemma~\ref{lem:csplift}, each
$(\widehat\omega_i,\widehat\eta_i)$ is already a symplectic pair;
if the original pair has type $(p_i,q_i)$, its lift has type
$(p_i+1,q_i)$. Thus the additional requirement of Definition \ref{def:csp} is  that
\begin{equation}\label{eq:cspforms}
  \widehat\Omega_i
  =\widehat\omega_i+\widehat\eta_i
  =d\!\left(\frac{r^2}{2}\pi^*\alpha_i\right)+\pi^*\omega_i, \qquad i=1,2,3
\end{equation}
form a positive quaternionic triple.
Let $(\widehat g,I_1,I_2,I_3)$ be the canonical hyper-K\"ahler
structure of the lifted hyper symplectic pair, provided by
Proposition~\ref{prop:sphyperkahler}, and put $R=r\partial_r$.
The hyper contact--symplectic pair is called
\emph{radially invariant} if $\Lie_RI_i=0$ for all $i$.
A radially invariant hyper pair is \emph{normalized} if
\begin{equation}\label{eq:cspnormalization}
  \widehat g(R,\cdot)=r\,dr.
\end{equation}
Its \emph{canonical metric} is the metric on $M$ induced by
the inclusion $M\cong M\times\{1\}\subset\widehat M$.

\begin{proposition}
\label{prop:cspspin}
Let $M$  be endowed with a hyper contact--symplectic pair. Then
  $\dim M\equiv 3\pmod 4$,
and   $TM\oplus\mathbb R$ carries a quaternionic Hermitian structure, where, with a slight
abuse of notation, $\mathbb R$ denotes the trivial real line
bundle over $M$.
\end{proposition}

\begin{proof}
By Definition~\ref{def:csp}, the lifted pairs on
$\widehat M=M\times\R_{>0}$ form a hyper symplectic pair.
Proposition~\ref{prop:sphyperkahler} gives a canonical hyper-K\"ahler
structure $(\widehat g,I_1,I_2,I_3)$ on $\widehat M$.
In particular, for cyclic $(i,j,k)$,
\begin{equation}\label{eq:quaternionic-identities-spin}
  I_i^2=-\operatorname{Id},
  \qquad I_iI_j=I_k=-I_jI_i,
\end{equation}
and
\begin{equation}\label{eq:quaternionic-hermitian-spin}
  \widehat g(I_iU,I_iV)=\widehat g(U,V).
\end{equation}
Thus $\dim\widehat M=4n$ for an integer $n$, and
 $ \dim M=\dim\widehat M-1=4n-1\equiv3\pmod4$.
Let $j_1:M\longrightarrow\widehat M$ be the inclusion
$j_1(x)=(x,1)$. Restricting the  ambient tangent bundle
 gives
  $j_1^*T\widehat M
  \cong TM\oplus\mathbb R\partial_r
  \cong TM\oplus\mathbb R.$
The restrictions of $\widehat g,I_1,I_2,I_3$ to this bundle
satisfy the same identities and provide its quaternionic
Hermitian structure.
\end{proof}

\begin{theorem}\label{thm:cspalignment}
A hyper contact--symplectic pair is radially invariant if and only if
\begin{equation}\label{eq:cspalignment}
 \charac(\alpha_1)=\charac(\alpha_2)=\charac(\alpha_3)=\F,
 \qquad
 \ker\omega_1=\ker\omega_2=\ker\omega_3=\E.
\end{equation}
In this case all three types are $(2l+1,2m)$ for some $l,m\ge0$, and
\[
 TM=\E\oplus\F,\qquad \dim\E=4l+3,\qquad\dim\F=4m.
\]
Locally $M=P\times N$, the $\alpha_i$ come from $P$, and the $\omega_i$ come from $N$. The lifted hyper-K\"ahler metric is a local product of a radially invariant hyper-K\"ahler symplectization of $P$ and a hyper-K\"ahler metric on $N$.
\end{theorem}
\begin{proof}
Let $(p_i,q_i)$ be the type of $(\alpha_i,\omega_i)$ and consider
the lift on $\widehat M=M\times\R_{>0}$
\[
  \widehat\omega_i=d\!\left(\frac{r^2}{2}\alpha_i\right),
  \qquad \widehat\eta_i=\omega_i,
  \qquad \widehat\Omega_i=\widehat\omega_i+\widehat\eta_i.
\]
By Definition~\ref{def:csp},
these pairs form a hyper symplectic pair. We write
$(\widehat g,I_1,I_2,I_3)$ for its canonical structure and
$R=r\partial_r$. Lemma~\ref{lem:csplift} gives, at every $(x,r)\in\widehat M$,
\begin{equation}\label{kernels}
  (\ker\widehat\omega_i)_{(x,r)}
  =\charac(\alpha_i)_x\oplus\{0\},
  \qquad
  (\ker\widehat\eta_i)_{(x,r)}
  =(\ker\omega_i)_x\oplus\mathbb R\partial_r.
\end{equation}
These kernels are complementary, and the ranks of the forms are
$2p_i+2$ and $2q_i$, respectively.
Since the original forms are pulled back from $M$,
$\Lie_R\alpha_i=\Lie_R\omega_i=0$. Lie differentiation commutes
with exterior differentiation, so
\[
  \Lie_R\widehat\omega_i
  =d\!\left(\Lie_R\left(\frac{r^2}{2}\alpha_i\right)\right)
  =2\widehat\omega_i,
  \qquad
  \Lie_R\widehat\eta_i=0.
\]
Theorem~\ref{thm:radial}, applied with
$A_i=\widehat\omega_i$ and $B_i=\widehat\eta_i$, therefore gives
\[
  \Lie_RI_i=0\quad(i=1,2,3)
  \quad\Longleftrightarrow\quad
  \widehat\Pi_1=\widehat\Pi_2=\widehat\Pi_3,
  \quad
  \widehat\Pi_i
  :=(\widehat\Omega_i^\flat)^{-1}\circ\widehat\omega_i^\flat.
\]
By Theorem~\ref{thm:spalignment}, equality of these projectors is
precisely the condition that the lifted hyper symplectic pair
have common characteristic distributions. Hence, in view of \eqref{kernels}, we get \eqref{eq:cspalignment}. This proves the required equivalence. Assume henceforth that these conditions hold and put
  $\E=\ker\omega_i$,
  $\F=\charac(\alpha_i)$, $TM=\E\oplus\F$.
The common distributions on the lift are
\[
  \widehat\E=\ker\widehat\eta_i
  =\pi^*\E\oplus\mathbb R\partial_r,
  \qquad
  \widehat\F=\ker\widehat\omega_i=\pi^*\F,
\]
where $\pi:\widehat M\longrightarrow M$ is the projection and
pullback subbundles are identified with their horizontal lifts.
Theorem~\ref{thm:spalignment} implies that these distributions
are orthogonal, quaternionic and parallel for $\widehat g$.
The rank--nullity theorem, using
$\dim\widehat M=2p_i+2q_i+2$, gives
\begin{align*}
  \operatorname{rank}\widehat\E
  =\dim\widehat M-\operatorname{rank}\widehat\eta_i
    =2p_i+2,\quad
  \operatorname{rank}\widehat\F
  =\dim\widehat M-\operatorname{rank}\widehat\omega_i
    =2q_i.
\end{align*}
The ranks are independent of $i$ and are divisible by $4$ by
Theorem~\ref{thm:spalignment}. Since $\widehat\E$ contains
$\partial_r$, its rank is positive. We may therefore write
  $2p_i+2=4(l+1)$, $2q_i=4m$
with $l,m\geq0$. Hence every original pair has type
$(2l+1,2m)$, and we have that
  $\operatorname{rank}\E=4l+3$
  and $\operatorname{rank}\F=4m$.
We describe the factors on $M$ in order to identify the product
on the lift stated by Theorem~\ref{thm:spalignment}.
The distribution $\E=\ker\omega_1$ is integrable because
$\omega_1$ is closed. For $X,Y\in\Gamma(\F)$,
  $\alpha_1([X,Y])=-d\alpha_1(X,Y)=0$.
Also, $\Lie_Xd\alpha_1=0$ by Cartan's formula, so  that
\[
  \iota_{[X,Y]}d\alpha_1
  =\Lie_X(\iota_Yd\alpha_1)-\iota_Y(\Lie_Xd\alpha_1)=0.
\]
This proves the integrability of $\F$.
As in the proof of Theorem~\ref{thm:spalignment}, Frobenius theorem gives
local submersions with kernels $\F$ and $\E$. Their product map
has invertible differential because $\E\cap\F=\{0\}$.
After restriction to a smaller neighborhood, we obtain
\[
  U\cong P^{4l+3}\times N^{4m},
  \qquad \E|_U=TP,\qquad \F|_U=TN.
\]
We take the two local factors connected.
For any $Y\in\Gamma(\F)$ and each $i$, we have
  $\alpha_i(Y)=0$, $\iota_Yd\alpha_i=0$, so that
$\Lie_Y\alpha_i=d(\alpha_i(Y))+\iota_Yd\alpha_i=0$.
These identities show that $\alpha_i$ has no component on the
$N$-factor and is independent of that factor. Explicitly, if
$X$ is lifted from $P$ and $Y$ from $N$, then
  $Y(\alpha_i(X))=(\Lie_Y\alpha_i)(X)+\alpha_i([Y,X])=0$.
Thus there is a contact form $\alpha_i^P$ on $P$ such that
$\alpha_i=\pr_P^*\alpha_i^P$. Similarly, the fact that, for any $X\in\Gamma(\E)$,
  $\iota_X\omega_i=0$ and $\Lie_X\omega_i=0$
show, by the pullback argument in Theorem~\ref{thm:spalignment},
that $\omega_i=\pr_N^*\omega_i^N$ for a closed $2$-form on $N$.
The nonvanishing of
$\alpha_i\wedge(d\alpha_i)^{2l+1}\wedge\omega_i^{2m}$
implies the nonvanishing of both
$\alpha_i^P\wedge(d\alpha_i^P)^{2l+1}$ on $P$ and
$(\omega_i^N)^{2m}$ on $N$. Thus the first form is indeed contact
and the second is symplectic. 
Under the induced identification
$U\times\R_{>0}\cong(P\times\R_{>0})\times N$,
the common characteristic distributions are
$\widehat\E=T(P\times\R_{>0})$ and $\widehat\F=TN$.
Moreover, the two families of lifted forms are pulled back from
these very factors and are given there by
$\Omega_i^P=d\!\left(\frac{r^2}{2}\alpha_i^P\right)$ and 
  $\omega_i^N$.
Theorem~\ref{thm:spalignment} now gives a hyper-K\"ahler structure
$(\widehat g_P,J_1,J_2,J_3)$ on $P\times\R_{>0}$, with fundamental
forms $\Omega_i^P$. It also gives a hyper-K\"ahler structure
$(g_N,K_1,K_2,K_3)$ on $N$, with fundamental forms $\omega_i^N$.
Their metrics and
endomorphisms are exactly the restrictions of the canonical
reconstruction, so
\[
  \widehat g=\widehat g_P+g_N,\qquad I_i=J_i\oplus K_i.
\]
Finally, on the first factor, $\Lie_R\Omega_i^P=2\Omega_i^P$.
To verify the radial invariance of the reconstructed
endomorphisms explicitly, write
  $P_+:=P\times\mathbb R_{>0}$
and, for $a\in\{1,2,3\}$, set
\[
  b_a:=(\Omega_a^P)^\flat:TP_+\longrightarrow T^*P_+,
  \qquad
  b_a(X)=\iota_X\Omega_a^P.
\]
Each $b_a$ is a smooth vector bundle isomorphism,
since $\Omega_a^P$ is nondegenerate.
For every smooth vector field $X$ on $P_+$, the
Lie derivative of $b_a$ is given by
  $(\mathcal L_Rb_a)(X)
  =
  \mathcal L_R\bigl(b_a(X)\bigr)-b_a([R,X])$.
Evaluating this $1$-form on an arbitrary smooth
vector field $Y$, we obtain
\begin{align*}
  \bigl((\mathcal L_Rb_a)(X)\bigr)(Y)
  &=
  R\bigl(\Omega_a^P(X,Y)\bigr)
  -\Omega_a^P(X,[R,Y])
  -\Omega_a^P([R,X],Y)\\
  &=
  (\mathcal L_R\Omega_a^P)(X,Y)\\
  &=
  2\Omega_a^P(X,Y)\\
  &=
  \bigl(2b_a(X)\bigr)(Y).
\end{align*}
Since $X,Y$ are arbitrary, this proves
  $\mathcal L_Rb_a=2b_a$.
We next compute the Lie derivative of the inverse
bundle map
$b_a^{-1}:T^*P_+\longrightarrow TP_+$.
Here, for any smooth $1$-form $\lambda$,
\[
  \bigl(\mathcal L_R(b_a^{-1})\bigr)(\lambda)
  =
  [R,b_a^{-1}(\lambda)]
  -b_a^{-1}(\mathcal L_R\lambda).
\]
Hence  $b_a^{-1}\circ b_a=\operatorname{Id}_{TP_+}$
gives
\begin{align*}
  0=
  \mathcal L_R
  \bigl(b_a^{-1}\circ b_a\bigr)=
  \mathcal L_R(b_a^{-1})\circ b_a
  +b_a^{-1}\circ(\mathcal L_Rb_a).
\end{align*}
The left-hand side is zero because
  $(\mathcal L_R\operatorname{Id}_{TP_+})(X)
  =
  [R,X]-[R,X]=0$.
Composing the preceding identity on the right
with $b_a^{-1}$, we find
\begin{align*}
  \mathcal L_R(b_a^{-1})
  =
  -b_a^{-1}\circ(\mathcal L_Rb_a)\circ b_a^{-1}=
  -b_a^{-1}\circ(2b_a)\circ b_a^{-1}=
  -2b_a^{-1}.
\end{align*}
Finally, for every cyclic permutation $(i,j,k)$,
the reconstruction formula is
  $J_i=-b_j^{-1}\circ b_k$.
A further application of the Leibniz rule yields
\begin{align*}
  \mathcal L_RJ_i
  =
  -\mathcal L_R(b_j^{-1})\circ b_k
  -b_j^{-1}\circ(\mathcal L_Rb_k)=
  2b_j^{-1}\circ b_k
  -2b_j^{-1}\circ b_k
  =0.
\end{align*}
Thus the reconstructed endomorphisms $J_i$ are
radially invariant. The metric formula
$\widehat g_P=-\Omega_i^P(\,\cdot\,,J_i\,\cdot\,)$ then gives
$\Lie_R\widehat g_P=2\widehat g_P$.
Thus the first factor is a radially invariant hyper-K\"ahler
symplectization of $P$, and the second factor is hyper-K\"ahler.
\end{proof}

\begin{lemma}\label{lem:cone}
Suppose the forms $d(r^2\alpha_i/2)$ define a positive quaternionic triple on $P\times\R_{>0}$. If its metric satisfies $\widehat g(R,\cdot)=r\,dr$, then
\begin{equation}\label{formula_cono}
 \widehat g=dr^2+r^2g_P
\end{equation}
and the $\alpha_i$ are the contact forms of a $3$-Sasakian structure on $(P,g_P)$.
\end{lemma}
\begin{proof}
Let
  $R=r\partial_r$
and denote by
 $\delta_t:P\times\mathbb R_{>0}\longrightarrow
  P\times\mathbb R_{>0}$,
  $\delta_t(x,r)=(x,e^t r)$,
the flow generated by $R$.
For each $i\in\{1,2,3\}$,
  $\widehat\Omega_i
  =
  d\!\left(\frac{r^2}{2}\alpha_i\right)$.
Since the forms $\alpha_i$ are pulled back from $P$, we have
$\delta_t^*\alpha_i=\alpha_i$
and therefore
  $\delta_t^*\widehat\Omega_i=
  d\!\left(
     \frac{e^{2t}r^2}{2}\alpha_i
  \right)=
  e^{2t}\widehat\Omega_i$.
Equivalently,
\begin{equation}\label{flusso1}
  \mathcal L_R\widehat\Omega_i
  =2\widehat\Omega_i.
\end{equation}
We next consider the almost complex structures reconstructed
algebraically from the triple. For every cyclic permutation
$(i,j,k)$,
  $I_i
  =
  -(\widehat\Omega_j^\flat)^{-1}
   \circ\widehat\Omega_k^\flat$.
Since all three $2$-forms are multiplied by the same factor
$e^{2t}$ under $\delta_t$, this factor cancels in the above
expression. Hence
  $\delta_t^*I_i=I_i$,
or, equivalently,
\begin{equation}\label{flusso2}
  \mathcal L_RI_i=0.
\end{equation}
The reconstructed metric is characterized by
  $\widehat g(U,V)
  =
  -\widehat\Omega_i(U,I_iV)$.
Using \eqref{flusso1} and \eqref{flusso2}, we obtain
  $\delta_t^*\widehat g=e^{2t}\widehat g$,
and consequently
\begin{equation}\label{lie_derivative_g}
  \mathcal L_R\widehat g=2\widehat g.
\end{equation}
We now prove \eqref{formula_cono}.
First, as $R=r\partial_r$ and $r>0$, the 
condition
  $\widehat g(R,\cdot)=r\,dr$
gives
  $\widehat g(\partial_r,\cdot)=dr$.
Consequently,
\begin{equation}\label{metrica_cono1}
  \widehat g(\partial_r,\partial_r)=1,
  \qquad
  \widehat g(\partial_r,U)=0
  \quad\text{for every }U\in\ker( dr).
\end{equation}
Thus the radial direction is unit and orthogonal to
the tangent spaces of the level sets of $r$.
Let $X,Y$ be arbitrary local vector fields on $P$,
and denote by $\widetilde X,\widetilde Y$ their lifts
to $P\times\mathbb R_{>0}$, defined by
  $\widetilde X_{(p,r)}=(X_p,0)$
  and
  $\widetilde Y_{(p,r)}=(Y_p,0)$
under the canonical product decomposition of the
tangent bundle.
These vector fields have zero radial component and are
independent of $r$. Hence
  $[\partial_r,\widetilde X]=0$ and
  $\widetilde X(r)=0$,
and therefore
\[
  [R,\widetilde X]
  =[r\partial_r,\widetilde X]
  =r[\partial_r,\widetilde X]
   -\widetilde X(r)\partial_r
  =0.
\]
The same holds for $\widetilde Y$.
Using \eqref{lie_derivative_g}, we obtain
\begin{align*}
  2\widehat g(\widetilde X,\widetilde Y)
  &=
  (\mathcal L_R\widehat g)
  (\widetilde X,\widetilde Y)\\
  &=
  R\bigl(\widehat g(\widetilde X,\widetilde Y)\bigr)
  -\widehat g([R,\widetilde X],\widetilde Y)
  -\widehat g(\widetilde X,[R,\widetilde Y])\\
  &=
  r\,\partial_r
  \bigl(\widehat g(\widetilde X,\widetilde Y)\bigr).
\end{align*}
As $r>0$, the preceding identity gives
  $\partial_r\bigl(\widehat g(\widetilde X,\widetilde Y)\bigr)
  =\frac{2}{r}\widehat g(\widetilde X,\widetilde Y)$.
Hence, by the product rule,
\begin{equation}\label{derivata1}
  \partial_r\!
  \left(r^{-2}\widehat g(\widetilde X,\widetilde Y)\right) 
=
  -2r^{-3}\widehat g(\widetilde X,\widetilde Y)
  +r^{-2}\partial_r
   \bigl(\widehat g(\widetilde X,\widetilde Y)\bigr)=
  -2r^{-3}\widehat g(\widetilde X,\widetilde Y)
  +r^{-2}\frac{2}{r}
   \widehat g(\widetilde X,\widetilde Y)=0.
\end{equation}
Define a Riemannian metric $g_P$ on $P$ by restricting
$\widehat g$ to the level set $r=1$:
\[
  (g_P)_p(v,w)
  :=
  \widehat g_{(p,1)}\bigl((v,0),(w,0)\bigr),
  \qquad v,w\in T_pP.
\]
By \eqref{derivata1} we have that, for any
fixed $p$, the function
\[
  r\longmapsto
  r^{-2}
  \widehat g_{(p,r)}
  \bigl(\widetilde X_{(p,r)},\widetilde Y_{(p,r)}\bigr)
\]
is constant on $\mathbb R_{>0}$.
Evaluating it at $r=1$, we conclude that
\begin{equation}\label{metrica_cono2}
  \widehat g_{(p,r)}
  \bigl(\widetilde X_{(p,r)},\widetilde Y_{(p,r)}\bigr)
  =r^2(g_P)_p(X_p,Y_p).
\end{equation}
Combining \eqref{metrica_cono2} with \eqref{metrica_cono1}, we obtain \eqref{formula_cono}.
Now we prove the final part of the Lemma. Define
  $\widehat\xi_i:=-I_iR$
on $P\times\mathbb R_{>0}$.
Since $I_i$ is skew-adjoint with respect to $\widehat g$,
  $\widehat g(R,\widehat\xi_i)
  =-\widehat g(R,I_iR)=0$.
Thus $\widehat\xi_i$ is tangent to the level sets of $r$.
Moreover, using \eqref{flusso2}, we have
  $[R,\widehat\xi_i]
  =-(\mathcal L_RI_i)R-I_i[R,R]=0$.
Since $\widehat\xi_i(r)=0$, this gives
  $0=[r\partial_r,\widehat\xi_i]
   =r[\partial_r,\widehat\xi_i]$.
Therefore $\widehat\xi_i$ is independent of $r$ and
is the lift of a vector field $\xi_i$ on $P$.
From now on, we use the same notation for vector fields
on $P$ and their $r$-independent lifts. In particular,
  $I_iR=-\xi_i$.
For a vector field $X$ lifted from $P$, the identity
\[
  \widehat\Omega_i
  =r\,dr\wedge\alpha_i+\frac{r^2}{2}d\alpha_i
\]
implies
\begin{align*}
  r^2\alpha_i(X)=
  \widehat\Omega_i(R,X)=
  \widehat g(R,I_iX)=
  -\widehat g(I_iR,X)=
  r^2g_P(\xi_i,X).
\end{align*}
It follows that
  $\alpha_i=g_P(\xi_i,\cdot)$.
The quaternionic identities and the orthogonality of
the $I_i$ also imply
  $\widehat g(I_iR,I_jR)=r^2\delta_{ij}$.
Consequently,
\[
  g_P(\xi_i,\xi_j)=\delta_{ij},
  \qquad
  \alpha_i(\xi_j)=\delta_{ij}.
\]
The radial component of $I_iX$ is determined by
\[
  \widehat g(R,I_iX)=r^2\alpha_i(X).
\]
Since $\widehat g(R,R)=r^2$, the vector field
  $I_iX-\alpha_i(X)R$
is tangent to the level sets of $r$.
It commutes with $R$, because
$\mathcal L_RI_i=0$, $[R,X]=0$, and
$R(\alpha_i(X))=0$. As above, it is therefore independent
of $r$. In this way we obtain a tensor field
$\phi_i:TP\longrightarrow TP$ characterized by
\[
  I_iX=\phi_iX+\alpha_i(X)R.
\]
Applying $I_i^2=-\operatorname{Id}$ to $X$ and comparing
tangential and radial components yields
\[
  \phi_i^2=-\operatorname{Id}_{TP}
           +\alpha_i\otimes\xi_i,
  \qquad
  \alpha_i\circ\phi_i=0.
\]
Likewise, $I_i\xi_i=R$ gives $\phi_i\xi_i=0$.
The compatibility of $\widehat g$ with $I_i$ implies
\[
  g_P(\phi_iX,\phi_iY)
  =g_P(X,Y)-\alpha_i(X)\alpha_i(Y).
\]
Thus $(\phi_i,\xi_i,\alpha_i,g_P)$ is an almost contact
metric structure.
For every cyclic permutation $(i,j,k)$, applying
$I_iI_j=I_k=-I_jI_i$ to $R$ and to tangent fields gives
\begin{gather*}
  \phi_k
  =\phi_i\phi_j-\alpha_j\otimes\xi_i
  =-\phi_j\phi_i+\alpha_i\otimes\xi_j, \\
  \xi_k=\phi_i\xi_j=-\phi_j\xi_i, \quad
  \alpha_k
  =\alpha_i\circ\phi_j
  =-\alpha_j\circ\phi_i.
\end{gather*}
Hence the three structures form an almost contact
metric $3$-structure.
Furthermore, for every vector fields $X,Y$ lifted from $P$,
\begin{align*}
  \frac{r^2}{2}d\alpha_i(X,Y)
  =
  \widehat\Omega_i(X,Y)=
  \widehat g(X,I_iY)=
  r^2g_P(X,\phi_iY).
\end{align*}
Therefore
\[
  d\alpha_i(X,Y)=2g_P(X,\phi_iY)
\]
and
$(\phi_i,\xi_i,\alpha_i,g_P)_{i=1,2,3}$ is a contact metric
$3$-structure. By Kashiwada's theorem \cite{Kashiwada01},
every contact metric $3$-structure is $3$-Sasakian.
Hence the $\alpha_i$ are the contact forms of a
$3$-Sasakian structure on $(P,g_P)$.
\end{proof}

\begin{theorem}\label{thm:cspnormalized}
A normalized hyper contact--symplectic pair $\left(\alpha_{i},\omega_{i}\right)_{i\in\left\{1,2,3\right\}}$ is locally, with its canonical metric, the product
\begin{equation}\label{eq:cspproduct}
 (M,g)=(P^{4l+3},g_P)\times(N^{4m},g_N),
\end{equation}
where $P$ is $3$-Sasakian and $N$ is hyper-K\"ahler. Conversely every such local product, with the forms globally defined, is a normalized hyper contact--symplectic pair.
In particular $M$ has a canonical $3$-quasi-Sasakian structure of rank $4l+3$. Its fundamental forms are
\begin{equation}\label{eq:qsfundamental}
 \Phi_i=\frac12d\alpha_i+\omega_i.
\end{equation}
Writing $Z_i$ for its Reeb fields, one has, for cyclic $(i,j,k)$,
\begin{equation}\label{eq:correctcontractions}
 \iota_{Z_j}\omega_i=0\quad\text{for all }i,j,
 \qquad
 \alpha_i=\frac12\iota_{Z_j}d\alpha_k
          =-\frac12\iota_{Z_k}d\alpha_j,
 \qquad [Z_i,Z_j]=2Z_k.
\end{equation}
\end{theorem}
\begin{proof}
Let
  $\widehat M=M\times\mathbb R_{>0}$ and
  $R=r\partial_r$,
and denote by $\widehat g$ the canonical metric of the lifted
hyper symplectic pair. Its fundamental forms are
\[
  \widehat\Omega_i
  =
  d\!\left(\frac{r^2}{2}\alpha_i\right)+\omega_i.
\]
Let $g$ be the canonical metric on $M$, obtained by
restricting $\widehat g$ to the hypersurface $r=1$.
By definition, the given hyper contact--symplectic pair
is radially invariant and satisfies
  $\widehat g(R,\cdot)=r\,dr$.
Thus, we can use
Theorem~\ref{thm:cspalignment}.
In particular, the distributions
  $\E=\ker\omega_i$
  and
  $\F=\ker\alpha_i\cap\ker d\alpha_i$
are independent of $i\in\{1,2,3\}$, and
\[
  TM=\E\oplus\F,
  \qquad
  \operatorname{rank}\E=4l+3,
  \qquad
  \operatorname{rank}\F=4m.
\]
Every point of $M$ has a neighborhood identified with
a product
  $U\cong P\times N$,
where $\E$ and $\F$ correspond to $TP$ and $TN$,
respectively. Under this identification,
\[
  \alpha_i=\alpha_i^P,
  \qquad
  \omega_i=\omega_i^N,
\]
where $\alpha_i^P$ is a contact form on $P$ and
$\omega_i^N$ is a two-form on $N$.
Theorem~\ref{thm:spalignment}, applied to the lifted hyper
symplectic pair as in the proof of
Theorem~\ref{thm:cspalignment}, gives a local Riemannian product
on the lift:
\[
  U\times\mathbb R_{>0}
  \cong
  (P\times\mathbb R_{>0})\times N,
  \qquad
  \widehat g=\widehat g_P+g_N.
\]
Here $\widehat g_P$ is the hyper-K\"ahler metric
reconstructed from
  $d\!\left(\frac{r^2}{2}\alpha_i^P\right)$,
whereas $g_N$ is the hyper-K\"ahler metric whose
fundamental forms are $\omega_i^N$.
Since $R$ is tangent to $P\times\mathbb R_{>0}$,
the normalization of $\widehat g$ restricts to
$\widehat g_P(R,\cdot)=r\,dr$.
Lemma~\ref{lem:cone} now implies
  $\widehat g_P=dr\otimes dr+r^2g_P$,
where $(P,g_P)$ is $3$-Sasakian and the $\alpha_i^P$
are its defining contact forms.
Consequently,
\[
  \widehat g=dr\otimes dr+r^2g_P+g_N.
\]
Restriction to $r=1$ gives
\[
  g=g_P+g_N.
\]
This proves the local Riemannian product
\eqref{eq:cspproduct}, with $P$ $3$-Sasakian and
$N$ hyper-K\"ahler.

Now, let
$(\phi_i^P,\xi_i^P,\alpha_i^P,g_P)$ denote the $3$-Sasakian structure on $P$, and
$(I_i^N,g_N)$ the hyper-K\"ahler structure on $N$, with
$\omega_i^N(U,V)=g_N(U,I_i^NV)$.
Define on $M$ the $2$-forms
\[
  \Phi_i:=\frac12d\alpha_i+\omega_i.
\]
Since $g$ is nondegenerate, there exist unique smooth
vector fields $Z_i$ and $(1,1)$ tensor fields
$\phi_i$ on $M$ such that
\[
  g(Z_i,X)=\alpha_i(X),
  \qquad
  g(X,\phi_iY)=\Phi_i(X,Y)
\]
for all vector fields $X,Y$ on $M$.
On each local product $P\times N$, the identities
\[
  d\alpha_i^P(X,Y)=2g_P(X,\phi_i^PY),
  \qquad
  \omega_i^N(U,V)=g_N(U,I_i^NV)
\]
show that these globally defined tensors have the
expressions
\[
  Z_i=(\xi_i^P,0),
  \qquad
  \phi_i=\phi_i^P\oplus I_i^N.
\]
The almost contact metric identities on $P$ and the
quaternionic Hermitian identities on $N$ therefore give
\begin{gather*}
  \phi_i^2=-\operatorname{Id}_{TM}
           +\alpha_i\otimes Z_i,
  \quad
  \alpha_i(Z_i)=1,\quad
  g(\phi_iX,\phi_iY)
  =g(X,Y)-\alpha_i(X)\alpha_i(Y).
\end{gather*}
For every cyclic permutation $(i,j,k)$, the same
decomposition implies
\begin{gather*}
  \phi_k
  =\phi_i\phi_j-\alpha_j\otimes Z_i
  =-\phi_j\phi_i+\alpha_i\otimes Z_j,
\\
  Z_k=\phi_iZ_j=-\phi_jZ_i,
  \quad
\alpha_k=\alpha_i\circ\phi_j=-\alpha_j\circ\phi_i.
\end{gather*}
On $TP$ these are the $3$-Sasakian identities,
while on $TN$ they reduce to the quaternionic
identities for the $I_i^N$.
Thus $(\phi_i,Z_i,\alpha_i,g)_{i=1,2,3}$ is an
almost contact metric $3$-structure on $M$.
Moreover,
\[
  \alpha_i(Z_i)=1,
  \qquad
  \iota_{Z_i}d\alpha_i=0,
  \qquad
  \iota_{Z_i}\omega_i=0,
\]
because locally $Z_i$ is the lift of $\xi_i^P$.
Hence $Z_i$ is also the Reeb vector field of the
original contact--symplectic pair $(\alpha_i,\omega_i)$.
Now we show normality of $(\phi_i,Z_i,\alpha_i)$ on $P\times N$, namely the vanishing of the tensor
\begin{align*}
  \mathcal N_i(X,Y):=
  [\phi_iX,\phi_iY]
  -\phi_i[\phi_iX,Y]
  -\phi_i[X,\phi_iY]
  +\phi_i^2[X,Y]
  +d\alpha_i(X,Y)Z_i.
\end{align*}
If $X,Y$ are lifted from $P$, the local expressions
for $\phi_i,Z_i,\alpha_i$ identify $\mathcal N_i(X,Y)$
with the lift of the normality tensor of
$(\phi_i^P,\xi_i^P,\alpha_i^P)$.
This vanishes because that structure is Sasakian.
If $X,Y$ are lifted from $N$, then
$d\alpha_i(X,Y)=0$ and $\phi_i=I_i^N$ on these fields.
Thus $\mathcal N_i(X,Y)$ is the lift of the
Nijenhuis tensor of $I_i^N$, which vanishes because
$(N,g_N,I_i^N)$ is K\"ahler.
Finally, suppose $X$ is lifted from $P$ and $Y$
is lifted from $N$. Then $\phi_iX$ is lifted from
$P$, while $\phi_iY$ is lifted from $N$.
Fields lifted from different factors commute, so
\[
  [X,Y]=[\phi_iX,Y]=[X,\phi_iY]
  =[\phi_iX,\phi_iY]=0.
\]
Also $d\alpha_i(X,Y)=0$, since $\alpha_i$ is pulled
back from $P$. Hence $\mathcal N_i(X,Y)=0$ in this
case as well.
The fundamental $2$-form is, by construction,
$\Phi_i=\frac12d\alpha_i+\omega_i$.
Since $d\omega_i=0$, we have
  $d\Phi_i=\frac12d^2\alpha_i+d\omega_i=0$.
Therefore each $(\phi_i,Z_i,\alpha_i,g)$ is
quasi-Sasakian, and  the
resulting almost contact $3$-structure is $3$-quasi-Sasakian.
To determine its rank, recall that $\alpha_i^P$ is
a contact form on the $(4l+3)$-dimensional factor $P$.
Consequently,
\[
  \alpha_i\wedge(d\alpha_i)^{2l+1}\neq0,
  \qquad
  (d\alpha_i)^{2l+2}=0.
\]
The first expression denotes a nowhere-zero
$(4l+3)$-form on $M$, obtained by pulling back the
contact volume form on $P$.
Hence each $\alpha_i$ has constant class $4l+3$,
which is the rank of the $3$-quasi-Sasakian structure.
Since any $Z_j$ is tangent to the $P$-factor, for all $i,j$ we have
  $\iota_{Z_j}\omega_i=0$.
Fix a cyclic permutation $(i,j,k)$.
The almost contact $3$-structure identities on $P$
give
  $\phi_k^P\xi_j^P=-\xi_i^P$
  and
  $\phi_j^P\xi_k^P=\xi_i^P$.
For an arbitrary vector field
$X=(X_P,X_N)$ on the local product, we therefore obtain
\begin{align*}
  d\alpha_k(Z_j,X)
=
  2g_P(\xi_j^P,\phi_k^PX_P)=
  -2g_P(\phi_k^P\xi_j^P,X_P)=
  2g_P(\xi_i^P,X_P)=2\alpha_i(X).
\end{align*}
Similarly,
  $d\alpha_j(Z_k,X)
  =
  -2g_P(\phi_j^P\xi_k^P,X_P)=
  -2\alpha_i(X)$.
Thus
\[
  \alpha_i
  =\frac12\iota_{Z_j}d\alpha_k
  =-\frac12\iota_{Z_k}d\alpha_j.
\]
Finally, for the Levi-Civita connection $\nabla^P$
of the Sasakian factor, we have
  $\nabla^P_X\xi_i^P=-\phi_i^PX$.
Hence
\begin{align*}
  [\xi_i^P,\xi_j^P]=
  \nabla^P_{\xi_i^P}\xi_j^P
  -\nabla^P_{\xi_j^P}\xi_i^P=
  -\phi_j^P\xi_i^P+\phi_i^P\xi_j^P=2\xi_k^P.
\end{align*}
Lifting this identity to $P\times N$ gives
  $[Z_i,Z_j]=2Z_k$.
All the identities in \eqref{eq:correctcontractions}
are therefore established.
Conversely, suppose that we have  a local Riemannian product
$(P^{4l+3},g_P)\times(N^{4m},g_N)$,
where $P$ is $3$-Sasakian and $N$ is hyper-K\"ahler.
Let $\alpha_i$ be the pullbacks of the defining
contact forms on $P$, and let $\omega_i$ be the
pullbacks of the hyper-K\"ahler forms on $N$.
For each $i\in\left\{1,2,3\right\}$,
\[
  \alpha_i\wedge(d\alpha_i)^{2l+1}
  \wedge\omega_i^{2m}
\]
is a volume form on the product, because it is the
exterior product of a contact volume form on $P$
and a symplectic volume form on $N$.
Furthermore,
\[
  (d\alpha_i)^{2l+2}=0,
  \qquad
  \omega_i^{2m+1}=0,
  \qquad
  d\omega_i=0.
\]
Hence $(\alpha_i,\omega_i)$ is a contact--symplectic
pair of type $(2l+1,2m)$.
The standard cone construction for a $3$-Sasakian
manifold (see \cite{BG99}) gives a hyper-K\"ahler structure on
$P\times\mathbb R_{>0}$ with metric
  $dr\otimes dr+r^2g_P$
and fundamental forms
  $d\!\left(\frac{r^2}{2}\alpha_i\right)$.
The converse part of Theorem~\ref{thm:spalignment}, applied to
this cone and the hyper-K\"ahler manifold $N$, shows that
\[
  \left(d\!\left(\frac{r^2}{2}\alpha_i\right),\omega_i\right)_{i=1,2,3}
\]
is a hyper symplectic pair. Its canonical metric is
\[
  \widehat g=dr\otimes dr+r^2g_P+g_N
\]
with fundamental forms
\[
  \widehat\Omega_i
  =
  d\!\left(\frac{r^2}{2}\alpha_i\right)+\omega_i.
\]
These are precisely the sums of the pairs prescribed in
Definition~\ref{def:csp}. Hence the original structures form a hyper
contact--symplectic pair.
The complex structures of the cone factor satisfy
$\mathcal L_RI_i^P=0$, as in the first part of the
proof of Lemma~\ref{lem:cone}.
The structures $I_i^N$ on the second factor are
independent of $r$. Therefore the product structures
satisfy
  $\mathcal L_RI_i=0$.
Moreover, recall that the metric constructed above is
$\widehat g=dr\otimes dr+r^2g_P+g_N$.
Here $g_P$ and $g_N$ denote the pullbacks of the metrics
on $P$ and $N$ to $(P\times\mathbb R_{>0})\times N$.
Since the projections of $R$ onto both $P$ and $N$
are zero, we have
  $g_P(R,U)=g_N(R,U)=0$
for every vector field $U$ on the product.
Consequently,
\begin{align*}
  \widehat g(R,U)=
  dr(R)\,dr(U)+r^2g_P(R,U)+g_N(R,U)=
  r\,dr(U).
\end{align*}
As $U$ is arbitrary, this proves that
  $\widehat g(R,\cdot)=r\,dr$.
Together with the radial invariance established above,
this shows that the hyper contact--symplectic pair
is normalized, and its canonical metric is
  $\widehat g|_{r=1}=g_P+g_N$.
 The algebraically reconstructed
metric and endomorphisms are uniquely determined
by these forms, so the local constructions agree
on overlaps. The positivity, radial invariance,
and normalization conditions hold on every such
neighborhood and hence hold globally.
This proves the converse.
\end{proof}

\begin{remark}\label{rem:3qs}
Every $3$-quasi-Sasakian structure gives rise to a hyper
contact--symplectic pair, but the construction differs
between the non-$3$-cosymplectic and the $3$-cosymplectic
cases. In particular, the resulting hyper pair need not
be normalized.
Let
$(M^{4n+3},\phi_i,\xi_i,\alpha_i,g_0)_{i=1,2,3}$
be a $3$-quasi-Sasakian manifold, and denote its
fundamental two-forms by
  $\Phi_i(X,Y):=g_0(X,\phi_iY)$.
The subscript in $g_0$ distinguishes the original metric
from the canonical metric of the hyper pair constructed
below. Recall that
\[
  [\xi_i,\xi_j]=c\,\xi_k
\]
for cyclic $(i,j,k)$ and a real constant $c$.

\medskip
\noindent
\textit{The case $c\neq0$.}
Put $\lambda=c/2$ and let $4l+3$ be the rank of the
$3$-quasi-Sasakian structure, with $m=n-l$.
By \cite[Corollary~5.7]{CD09}, $M$ is locally the
Riemannian product
\[
  (M,g_0)
  \cong
  (P^{4l+3},g_P)\times(N^{4m},g_N),
\]
where $P$ is $3$-$\lambda$-Sasakian and $N$ is
hyper-K\"ahler.
The original tensors respect this local decomposition.
In particular, $\alpha_i$ is pulled back from $P$ and,
with the exterior derivative convention adopted here,
\[
  d\alpha_i=2\lambda\Phi_i^P,
  \qquad
  \Phi_i=\Phi_i^P+\Phi_i^N.
\]
Here the forms on the two factors are understood as
pulled back to the product.
Define the global $2$-forms on $M$
\[
  \Theta_i
  :=
  \Phi_i-\frac{1}{2\lambda}d\alpha_i.
\]
These forms are closed, since $d\Phi_i=0$ and
$d^2\alpha_i=0$.
On every local product neighborhood,
  $\Theta_i=\Phi_i^N$.
Thus $\Theta_i$ is the hyper-K\"ahler form of the
$N$-factor, extended by zero on the $P$-factor.
In particular, the ordinary contact--symplectic pairs
  $\left\{(\alpha_i,\Theta_j)\right\}_{i,j=1,2,3}$
have type $(2l+1,2m)$. Note that these are the pairs already
described in \cite[Corollary~5.9]{CD09}.
To obtain a normalized hyper contact--symplectic pair,
set
  $\alpha_i':=\lambda\alpha_i$ and
  $\omega_i:=\Theta_i$.
Indeed, on the $P$-factor the rescaled tensors
\[
  \phi_i'=\phi_i^P,\qquad
  \xi_i'=\lambda^{-1}\xi_i^P,\qquad
  \alpha_i'=\lambda\alpha_i^P,\qquad
  g_P'=\lambda^2g_P
\]
form a $3$-Sasakian structure. In particular,
\[
  d\alpha_i'
  =2\lambda^2g_P(\,\cdot\,,\phi_i^P\,\cdot\,)
  =2g_P'(\,\cdot\,,\phi_i'\,\cdot\,).
\]
The converse part of Theorem~\ref{thm:cspnormalized}
therefore applies to
  $(P,g_P')\times(N,g_N)$.
It follows that
\[
  (\alpha_i',\omega_i)_{i=1,2,3}
  =
  (\lambda\alpha_i,\Theta_i)_{i=1,2,3}
\]
is a normalized hyper contact--symplectic pair of
type $(2l+1,2m)$.
More explicitly, the lifted forms
\[
  \widehat\Omega_i
  =
  d\!\left(\frac{r^2}{2}\alpha_i'\right)+\Theta_i
\]
are the hyper-K\"ahler forms of the local product
of the cone over $(P,g_P')$ with $(N,g_N)$.
Their reconstructed metric is
  $\widehat g
  =
  dr\otimes dr+r^2\lambda^2g_P+g_N$.
Hence the canonical metric on $M$ is locally
\[
  g_{\mathrm{can}}=\lambda^2g_P+g_N.
\]
This agrees with the original metric
$g_0=g_P+g_N$ precisely when $\lambda^2=1$.
In particular, when $c=2$, so that $\lambda=1$,
no rescaling is required:
  $(\alpha_i,\Theta_i)_{i=1,2,3}$
is already normalized and $g_{\mathrm{can}}=g_0$.

\medskip
\noindent
\textit{The case $c=0$.}
In this case the original structure is
$3$-cosymplectic (see \cite[Theorem~3.7]{CD09}).
Thus
  $d\alpha_i=0$ and  $d\Phi_i=0$.
Here we use the full fundamental forms $\Phi_i$,
rather than the forms $\Theta_i$ of the preceding
construction, and consider
\[
  (\alpha_i,\omega_i):=(\alpha_i,\Phi_i).
\]
Since each $(\phi_i,\xi_i,\alpha_i,g_0)$ is an almost
contact metric structure,
  $\alpha_i\wedge\Phi_i^{\,2n+1}$
is a volume form on $M$, while
  $\Phi_i^{\,2n+2}=0$.
Consequently, each $(\alpha_i,\Phi_i)$ is a
contact--symplectic pair of type $(0,2n+1)$.
We verify that these three pairs satisfy
Definition~\ref{def:csp}.
Their lifted forms are
\[
  \widehat\Omega_i
  =
  d\!\left(\frac{r^2}{2}\alpha_i\right)+\Phi_i
  =
  r\,dr\wedge\alpha_i+\Phi_i.
\]
On $M\times\mathbb R_{>0}$ define
\[
  \widehat g=r^2\,dr\otimes dr+g_0,
  \qquad
  T=\frac1r\partial_r,
\]
so that $T$ is unit and orthogonal to the
$M$-directions. Define endomorphisms by
\[
  I_iT=-\xi_i,
  \qquad
  I_iX=\phi_iX+\alpha_i(X)T
\]
for vector fields $X$ lifted from $M$.
The almost contact metric $3$-structure identities
give
\begin{equation*}
  I_i^2=-\operatorname{Id},
  \quad
  I_iI_j=I_k=-I_jI_i, \quad
  \widehat g(I_iU,I_iV)=\widehat g(U,V).
\end{equation*}
Moreover, if $U=X+aT$ and $V=Y+bT$, then
\begin{align*}
  \widehat g(U,I_iV)
  &=
  \Phi_i(X,Y)+a\alpha_i(Y)-b\alpha_i(X)\\
  &=
  (r\,dr\wedge\alpha_i+\Phi_i)(U,V)\\
  &=
  \widehat\Omega_i(U,V).
\end{align*}
Thus the lifted forms constitute a positive quaternionic
triple. By Lemma~\ref{lem:csplift}, they are the sums of the
symplectic pairs
\[
  \left(d\!\left(\frac{r^2}{2}\alpha_i\right),\Phi_i\right).
\]
These pairs therefore form a hyper symplectic pair, whose canonical
metric is precisely $\widehat g$. By
Proposition~\ref{prop:sphyperkahler}, the reconstructed structure
is hyper-K\"ahler. Applying Definition~\ref{def:csp}, we conclude that
  $(\alpha_i,\Phi_i)_{i=1,2,3}$
is a hyper contact--symplectic pair, and its
canonical metric is the original metric
  $g_{\mathrm{can}}=g_0$.
This hyper pair is not radially invariant.
Indeed, with $R=r\partial_r=r^2T$, we have
 $I_iR=-r^2\xi_i$.
The fields $\xi_i$ are independent of $r$, and hence
  $(\mathcal L_RI_i)(R)
  =
  [R,I_iR]-I_i[R,R]
  =
  -2r^2\xi_i\neq0$.
It is therefore not normalized. Its
metric also satisfies
  $\widehat g(R,\cdot)=r^3\,dr$,
\end{remark}

In particular, by Remark \ref{rem:3qs},  \cite{CDY13} and \cite{CDY15} provide  non-trivial examples of  hyper contact-symplectic manifolds.  We now construct compact examples for which the canonical local product decomposition described in Theorem~\ref{thm:cspnormalized}
does not extend to a global product.

Let $B$ be a connected manifold, let $\widetilde B\to B$ be its
universal covering, on which $\pi_1(B)$ acts by deck
transformations, and let $\rho:\pi_1(B)\to\operatorname{Diff}(F)$ be
a homomorphism. The \emph{foliated bundle} with holonomy $\rho$ is
\[
 \widetilde B\times_\rho F:=(\widetilde B\times F)/\pi_1(B),
 \qquad
 \gamma\cdot(b,y)=(\gamma b,\rho(\gamma)y).
\]
The action is free and properly discontinuous, thus $\widetilde B\times_\rho F$ is a manifold,
which fibres over $B$ with fibre $F$, and the images of the slices
$\widetilde B\times\{y\}$ are the leaves of a foliation transverse
to the fibres (see \cite{CN, C} for a more detailed description). It is compact if $B$ and $F$ are compact.

For a $3$-Sasakian manifold $(P,\phi_i,\xi_i,\alpha_i,g_P)$ we denote
by $\operatorname{Aut}(P)$ the group of isometries $f$ of $g_P$ such
that $f^*\alpha_i=\alpha_i$ for $i=1,2,3$; such maps also preserve the
$\xi_i$ and the $\phi_i$. For a hyper-K\"ahler manifold
$(N,g_N,\omega_i)$ we denote by $\operatorname{Aut}(N)$ the group of
isometries $f$ with $f^*\omega_i=\omega_i$ for $i=1,2,3$. 
For the standard $3$-Sasakian sphere $S^{4n+3}\subset\HH^{n+1}$, whose Reeb
fields are generated by right multiplication by $i,j,k$, the group
$\Sp(n+1)$, acting on the left, is contained in
$\operatorname{Aut}(S^{4n+3})$ (see \cite{BG99}). For a flat torus
$T^{4m}=\HH^m/\Lambda$ with the hyper-K\"ahler structure induced from
$\HH^m$, all translations and the involution $v\mapsto-v$ belong to
$\operatorname{Aut}(T^{4m})$, since the forms $\omega_i$ have
constant coefficients.

\begin{proposition}\label{prop:cspfoliated}
Let $(P^{4l+3},g_P)$ be a $3$-Sasakian manifold with contact forms
$\alpha_i$, and let $(N^{4m},g_N)$ be a hyper-K\"ahler manifold with
fundamental forms $\omega_i$. Let $\{B,F\}=\{P,N\}$ and let
$\rho:\pi_1(B)\to\operatorname{Aut}(F)$ be a homomorphism. Then the
pullbacks of the $\alpha_i$ and of the $\omega_i$ to
$\widetilde B\times F$ descend to forms, still denoted
$\alpha_i,\omega_i$, on $M=\widetilde B\times_\rho F$, and
$(\alpha_i,\omega_i)_{i=1,2,3}$ is a normalized hyper
contact--symplectic pair. Its canonical metric is induced by the
product metric, $\E=\ker\omega_i$ is tangent to the $P$-directions,
and $\F=\charac(\alpha_i)$ to the $N$-directions. If $P$ and $N$ are
compact, so is $M$.
\end{proposition}

\begin{proof}
Deck transformations preserve all the structure pulled back to
$\widetilde B$, and $\rho$ takes values in $\operatorname{Aut}(F)$.
Hence the forms and the product metric on $\widetilde B\times F$ are
invariant under the action of $\pi_1(B)$, and they descend to $M$.
Every point of $M$ has a neighborhood isometric to an open subset of
$\widetilde B\times F$, on which the descended forms are those of a
local product of a $3$-Sasakian and a hyper-K\"ahler manifold. The
converse part of Theorem~\ref{thm:cspnormalized} shows that
$(\alpha_i,\omega_i)$ is a normalized hyper contact--symplectic pair
with canonical metric induced by the product metric. The
description of $\E$ and $\F$ follows from
Theorem~\ref{thm:cspalignment}.
\end{proof}

\begin{example}\label{ex:cspfinite}
Let $l\geq0$ and $m\geq1$ be integers, and let
$\Gamma=\{\pm1\}\subset\Sp(1)$ act on $S^{4l+3}\subset\HH^{l+1}$ by
$x\mapsto-x$ and on a flat torus $T^{4m}$ by $v\mapsto-v$. Both maps
are automorphisms, and the first action is free. Taking
$B=\R P^{4l+3}=S^{4l+3}/\Gamma$, which is $3$-Sasakian, $F=T^{4m}$ and
$\rho$ the second action, Proposition~\ref{prop:cspfoliated} gives a
normalized hyper contact--symplectic pair on the compact manifold
\[
 M=(S^{4l+3}\times T^{4m})/\Z_2,\qquad (x,v)\sim(-x,-v),
\]
of dimension $4(l+m)+3$. The leaf of $\F$ through the image of
$(x,v)$ is the image of $\{x\}\times T^{4m}$, which is a torus
$T^{4m}$. The leaf of $\E$ through the image of $(x,v)$ is a copy of
$\R P^{4l+3}$ if $v$ is one of the $2^{4m}$ half-periods, which
satisfy $v=-v$, and a copy of $S^{4l+3}$ otherwise. Since leaves of
the same factor of a Riemannian product are all diffeomorphic, the
splitting $TM=\E\oplus\F$ is not a global product. More generally,
one may take $B=S^{4l+3}/\Gamma$ for any finite $\Gamma\subset\Sp(1)$
acting by left scalar multiplication, which is free, and any
homomorphism $\rho:\Gamma\to\operatorname{Aut}(N)$ to the automorphism
group of a compact hyper-K\"ahler manifold.
\end{example}

\begin{example}\label{ex:cspinfinite}
Let $l\geq0$ and $m\geq1$ be integers, let
$B=T^{4m}=\R^{4m}/\Z^{4m}$ with a flat hyper-K\"ahler structure
induced from $\HH^m$, let $F=S^{4l+3}$ be the standard $3$-Sasakian
sphere, and let $\rho:\Z^{4m}\to\Sp(l+1)$ be defined by
\[
 \rho(e_1)=e^{2\pi i\theta}\,\Id,\qquad \rho(e_a)=\Id\quad(a\geq2),
\]
with $e^{2\pi i\theta}\,\Id$ acting on $\HH^{l+1}$ by left
multiplication and $\theta\in\R\setminus\mathbb Q$. By
Proposition~\ref{prop:cspfoliated}, the compact manifold
$M=\R^{4m}\times_\rho S^{4l+3}$ carries a normalized hyper
contact--symplectic pair. The leaf of $\F=\charac(\alpha_i)$ through
$[v,y]$ is the image of $\R^{4m}\times\{y\}$, and the stabilizer of
$y$ in $\Z^{4m}$ is $\{n:n_1=0\}$, because $e^{2\pi in_1\theta}y=y$
forces $n_1=0$. Hence every leaf of $\F$ is diffeomorphic to
$\R\times T^{4m-1}$ and is not compact; its closure is the image of
$\R^{4m}\times\{e^{it}y:t\in\R\}$.
The canonical splitting of $TM$ does not become a global product after
passage to any finite covering of $M$. Indeed, the lifted $\F$-leaves remain noncompact,
whereas the factor leaves of a product of compact manifolds are compact.
\end{example}

\section{Hyper contact pairs}\label{sec:cp}

A \emph{contact pair} $(\alpha,\beta)$ of type $(p,q)$ on a manifold
$M$ of dimension $2p+2q+2$ consists of two $1$-forms such that
\begin{equation}\label{eq:cpordinary}
 \alpha\wedge(d\alpha)^p\wedge\beta\wedge(d\beta)^q\ne0,
 \qquad (d\alpha)^{p+1}=(d\beta)^{q+1}=0.
\end{equation}
In particular $\rank d\alpha=2p$ and $\rank d\beta=2q$. The
\emph{characteristic distributions} are
\[
 \charac(\alpha)=\ker\alpha\cap\ker d\alpha,
 \qquad
 \charac(\beta)=\ker\beta\cap\ker d\beta .
\]
The \emph{Reeb fields} $Z,W$ of the pair are uniquely determined by
\[
 \alpha(Z)=\beta(W)=1,\quad\alpha(W)=\beta(Z)=0,
 \quad \iota_Zd\alpha=\iota_Zd\beta
 =\iota_Wd\alpha=\iota_Wd\beta=0.
\]
Since $\ker d\alpha$ has rank $2q+2$ and contains $Z$, on which
$\alpha$ does not vanish, $\charac(\alpha)$ has rank $2q+1$;
similarly $\charac(\beta)$ has rank $2p+1$. Moreover
$W\in\charac(\alpha)$ and $Z\in\charac(\beta)$. If
$u\in\charac(\alpha)\cap\charac(\beta)$, then $u$ annihilates
$\alpha,\beta,d\alpha,d\beta$, so the contraction of $u$ with
the volume form \eqref{eq:cpordinary} vanishes and $u=0$. Since the
ranks add up to $\dim M$,
\begin{equation}\label{eq:cpsplitting}
 TM=\charac(\alpha)\oplus\charac(\beta).
\end{equation}
We refer to \cite{BH05,BH10} for more details.

\begin{lemma}\label{lem:cplift}
For a contact pair of type $(p,q)$, set on $\widehat M=M\times\R_{>0}^2$,
with coordinates $(r,s)$ on the second factor,
\begin{equation}\label{eq:cpAB}
 A=d\!\left(\frac{r^2}{2}\alpha\right),\qquad
 B=d\!\left(\frac{s^2}{2}\beta\right).
\end{equation}
Then
\begin{equation}\label{eq:cpkernels}
 \ker A=\charac(\alpha)\oplus\R\partial_s,
 \qquad\ker B=\charac(\beta)\oplus\R\partial_r.
\end{equation}
The kernels are complementary and the ranks are $2p+2$ and
$2q+2$, respectively. Thus $(A,B)$ is a symplectic pair of type
$(p+1,q+1)$, and $A+B$ is symplectic.
\end{lemma}

\begin{proof}
As in the proof of Lemma~\ref{lem:csplift}, we use the same
notation for forms on $M$ and for their pullbacks to $\widehat M$,
and we identify $T_{(x,r,s)}\widehat M$ with
$T_xM\oplus\R\partial_r\oplus\R\partial_s$. We have
\begin{equation}\label{eq:cpAB-explicit}
 A=r\,dr\wedge\alpha+\frac{r^2}{2}\,d\alpha,
 \qquad
 B=s\,ds\wedge\beta+\frac{s^2}{2}\,d\beta .
\end{equation}
Let $V=X+a\partial_r+b\partial_s$, with $X\in T_xM$ and
$a,b\in\R$. By \eqref{eq:cpAB-explicit},
\begin{equation}\label{prodotto-interno1}
 \iota_VA
  =ra\,\alpha-r\alpha(X)\,dr+\frac{r^2}{2}\iota_Xd\alpha .
\end{equation}
Assume $\iota_VA=0$. Evaluating on $\partial_r$ gives
$0=-r\alpha(X)$, hence $\alpha(X)=0$. Evaluating on the Reeb field
$Z$ and using $\alpha(Z)=1$, $dr(Z)=0$ and $\iota_Zd\alpha=0$ gives
$0=ra$, hence $a=0$. Then \eqref{prodotto-interno1} reduces to
$\frac{r^2}{2}\iota_Xd\alpha=0$, so $X\in\ker d\alpha$. Thus
$X\in\charac(\alpha)$ and $a=0$, while $b$ is arbitrary, because
$A$ has no $ds$-component. Conversely, every such vector lies in
$\ker A$ by \eqref{prodotto-interno1}. This proves
$\ker A=\charac(\alpha)\oplus\R\partial_s$.
In the same way, we obtain 
$\ker B=\charac(\beta)\oplus\R\partial_r$.
By \eqref{eq:cpsplitting},
\[
 \ker A\oplus\ker B
 =\charac(\alpha)\oplus\charac(\beta)\oplus\R\partial_r\oplus\R\partial_s
 =T\widehat M .
\]
Since $\rank\charac(\alpha)=2q+1$ and $\rank\charac(\beta)=2p+1$,
the kernels have ranks $2q+2$ and $2p+2$; as
$\dim\widehat M=2p+2q+4$, we get $\rank A=2p+2$ and
$\rank B=2q+2$.
Since $(dr\wedge\alpha)^2=0$ and $(d\alpha)^{p+1}=0$, the binomial
expansion of \eqref{eq:cpAB-explicit} gives, exactly as in
\eqref{eq:csp-A-power},
\[
 A^{p+1}=\frac{(p+1)\,r^{2p+1}}{2^p}\,dr\wedge\alpha\wedge(d\alpha)^p,
 \qquad A^{p+2}=0,
\]
and likewise
\[
 B^{q+1}=\frac{(q+1)\,s^{2q+1}}{2^q}\,ds\wedge\beta\wedge(d\beta)^q,
 \qquad B^{q+2}=0.
\]
Therefore
\[
 A^{p+1}\wedge B^{q+1}
 =-\frac{(p+1)(q+1)\,r^{2p+1}s^{2q+1}}{2^{p+q}}\,
  dr\wedge ds\wedge\alpha\wedge(d\alpha)^p\wedge\beta\wedge(d\beta)^q ,
\]
which is nowhere zero by \eqref{eq:cpordinary}. Both forms are
closed, being exact. Hence $(A,B)$ is a symplectic pair of type
$(p+1,q+1)$. Finally, in the expansion of $(A+B)^{p+q+2}$, as in \eqref{somma-simplettica}, only the
term containing $A^{p+1}\wedge B^{q+1}$ survives, 
and then $A+B$ is symplectic.
\end{proof}

\begin{definition}\label{def:cp}
Let $(\alpha_i,\beta_i)$, $i=1,2,3$, be three contact pairs on
$M$, and let
  $\pi:\widehat M=M\times\R_{>0}^2\longrightarrow M$
be the projection. The triple is called a \emph{hyper contact pair}
if the pairs
\[
  \widehat\omega_i
  :=d\!\left(\frac{r^2}{2}\pi^*\alpha_i\right),
  \qquad
  \widehat\eta_i
  :=d\!\left(\frac{s^2}{2}\pi^*\beta_i\right),
  \qquad i=1,2,3,
\]
form a hyper symplectic pair on $\widehat M$ in the sense of
Definition~\ref{def:sp}.
\end{definition}

By Lemma~\ref{lem:cplift}, each $(\widehat\omega_i,\widehat\eta_i)$
is already a symplectic pair. If the original pair has type
$(p_i,q_i)$, its lift has type $(p_i+1,q_i+1)$. Thus the
additional requirement of Definition~\ref{def:cp} is that
\begin{equation}\label{eq:cpforms}
  \widehat\Omega_i
  =\widehat\omega_i+\widehat\eta_i
  =d\!\left(\frac{r^2}{2}\pi^*\alpha_i
            +\frac{s^2}{2}\pi^*\beta_i\right), \qquad i=1,2,3,
\end{equation}
form a positive quaternionic triple.

Let $(\widehat g,I_1,I_2,I_3)$ be the canonical hyper-K\"ahler
structure of the lifted hyper symplectic pair, given by
Proposition~\ref{prop:sphyperkahler}, and put
$R=r\partial_r$, $S=s\partial_s$.
The hyper contact pair is called \emph{radially invariant} if
$\Lie_RI_i=0$ for every $i$, and \emph{normalized} if it is
radially invariant and
\begin{equation}\label{eq:cpnormalization}
  \widehat g(R,\cdot)=r\,dr,
  \qquad \widehat g(S,\cdot)=s\,ds.
\end{equation}
Its \emph{canonical metric} is the metric induced on the
submanifold $M\cong M\times\{(1,1)\}\subset\widehat M$.

\begin{remark}\label{rem:cpRS}
The definition of radial invariance does not privilege $R$ over
$S$: for a hyper contact pair, 
 $\Lie_RI_i=0$
 if and only if
 $\Lie_SI_i=0$ for all $i=1,2,3$.
Indeed, the flow of $R+S$ is $(x,r,s)\mapsto(x,e^tr,e^ts)$, which
multiplies all three forms \eqref{eq:cpforms} by the same factor
$e^{2t}$. Like at the beginning of the proof of
Lemma~\ref{lem:cone}, this factor cancels in
$I_i=-(\widehat\Omega_j^\flat)^{-1}\circ\widehat\Omega_k^\flat$, so
$\Lie_{R+S}I_i=0$, that is, $\Lie_SI_i=-\Lie_RI_i$. Equivalently, one
may apply Theorem~\ref{thm:radial} once with $(A_i,B_i,R)=
(\widehat\omega_i,\widehat\eta_i,R)$ and once with
$(A_i,B_i,R)=(\widehat\eta_i,\widehat\omega_i,S)$: both give the same
condition (iii).
\end{remark}

\begin{proposition}\label{prop:cpspin}
Let $(\alpha_i,\beta_i)_{i=1,2,3}$ be a hyper contact pair on $M$.
Then
 $\dim M=4n+2$
  for some $n\geq1$,
and $TM\oplus\R^2$ carries a quaternionic Hermitian structure,
where $\R^2$ denotes the trivial real plane bundle over $M$.
Moreover, the six $1$-forms
$\alpha_1,\alpha_2,\alpha_3,\beta_1,\beta_2,\beta_3$
are pointwise linearly independent.
\end{proposition}

\begin{proof}
By Proposition~\ref{prop:sphyperkahler}, the double symplectization
$\widehat M=M\times\R_{>0}^2$ carries its canonical hyper-K\"ahler
structure $(\widehat g,I_1,I_2,I_3)$. Hence
$\dim\widehat M\equiv0\pmod4$ and $\dim M\equiv2\pmod4$.
Restricting the tangent bundle and these tensors to the unit section
$j_{1,1}(x)=(x,1,1)$ gives
\[
 j_{1,1}^*T\widehat M
 \cong TM\oplus\R\partial_r\oplus\R\partial_s
 \cong TM\oplus\R^2,
\]
with the induced quaternionic Hermitian structure.
Along the unit section, the radial fields $R=r\partial_r$ and
$S=s\partial_s$ satisfy
\[
 \widehat\Omega_i(R,S)=0,\qquad
 \iota_R\widehat\Omega_i=\pi^*\alpha_i,\qquad
 \iota_S\widehat\Omega_i=\pi^*\beta_i.
\]
Thus
 $\pi^*\alpha_i=-\widehat g(I_iR,\cdot)$ and 
 $\pi^*\beta_i=-\widehat g(I_iS,\cdot)$.
The quaternionic Hermitian identities imply
\[
 \widehat g(I_iR,I_jR)=\delta_{ij}\widehat g(R,R),\qquad
 \widehat g(I_iS,I_jS)=\delta_{ij}\widehat g(S,S),
\]
and, using $\widehat\Omega_k(R,S)=0$,
\[
 \widehat g(I_iR,I_jS)=\delta_{ij}\widehat g(R,S).
\]
Consequently, the Gram matrix of the six vectors
$I_1R,I_1S,I_2R,I_2S,I_3R,I_3S$ consists of three copies of
\[
 \begin{pmatrix}
 \widehat g(R,R)&\widehat g(R,S)\\
 \widehat g(R,S)&\widehat g(S,S)
 \end{pmatrix}.
\]
This matrix is positive definite because $R$ and $S$ are linearly
independent. The six vectors, and hence their metric duals, are
therefore linearly independent. Since the latter are pullbacks
of the six forms on $M$, those forms are pointwise independent.
In particular $\dim M\geq6$, which completes the proof.
\end{proof}

\begin{remark}\label{rem:cp-lcs}
Suppose that the three contact pairs $(\alpha_i,\beta_i)$
have types $(h_i,0)$. Since $\dim M=2h_i+2$ for every $i$,
we have $h_1=h_2=h_3=:h$. By
Proposition~\ref{prop:cpspin}, $h=2n\geq2$.
Each pair determines the $2$-form
\[
 \sigma_i=d\alpha_i+\alpha_i\wedge\beta_i.
\]
Since $d\beta_i=0$, we have
\[
 d\sigma_i=\beta_i\wedge\sigma_i,\qquad
 \sigma_i^{h+1}
 =(h+1)\alpha_i\wedge(d\alpha_i)^h\wedge\beta_i\neq0.
\]
Thus $\sigma_i$ is locally conformally symplectic with Lee form
$\beta_i$, as in \cite[Proposition~3.1]{BK11}.
By Proposition~\ref{prop:cpspin}, the three Lee forms are
pointwise linearly independent. This gives a different
configuration from the locally conformal hyper-K\"ahler
structures in \cite{OP97}, whose fundamental forms satisfy
$d\Omega_i=\theta\wedge\Omega_i$ for a common closed Lee form
$\theta$.
\end{remark}

\begin{theorem}\label{thm:cpalignment}
A hyper contact pair is radially invariant if and only if
\begin{equation}\label{eq:cpalignment}
 \charac(\alpha_1)=\charac(\alpha_2)=\charac(\alpha_3)=\F,
 \qquad
 \charac(\beta_1)=\charac(\beta_2)=\charac(\beta_3)=\E.
\end{equation}
In this case all three pairs have the same type $(2l+1,2m+1)$ for
some $l,m\ge0$, and
\[
 TM=\E\oplus\F,\qquad \rank\E=4l+3,\qquad \rank\F=4m+3.
\]
Locally $M=P\times Q$, all the $\alpha_i$ are pulled back from
contact forms on $P$, and all the $\beta_i$ from contact forms on
$Q$. The lifted hyper-K\"ahler metric is the Riemannian product of
radially invariant hyper-K\"ahler symplectizations of $P$ and of $Q$.
\end{theorem}

\begin{proof}
Let $(p_i,q_i)$ be the type of the $i$th contact pair. By
Definition~\ref{def:cp}, the forms
\[
  \widehat\omega_i=d\!\left(\frac{r^2}{2}\alpha_i\right),
  \qquad
  \widehat\eta_i=d\!\left(\frac{s^2}{2}\beta_i\right)
\]
define a hyper symplectic pair on $\widehat M=M\times\R_{>0}^2$.
By Lemma~\ref{lem:cplift}, at every point of $\widehat M$,
\begin{equation}\label{eq:cpliftkernels}
  \ker\widehat\omega_i
  =\charac(\alpha_i)\oplus\R\partial_s,
  \qquad
  \ker\widehat\eta_i
  =\charac(\beta_i)\oplus\R\partial_r ;
\end{equation}
these kernels are complementary, and the forms have constant ranks
$2p_i+2$ and $2q_i+2$. Since $\alpha_i$, $\beta_i$ and $s$ do not
depend on $r$,
\[
  \Lie_R\widehat\omega_i
  =d\!\left(\Lie_R\left(\frac{r^2}{2}\alpha_i\right)\right)
  =2\widehat\omega_i,
  \qquad \Lie_R\widehat\eta_i=0.
\]
Hence Theorem~\ref{thm:radial} applies with $A_i=\widehat\omega_i$
and $B_i=\widehat\eta_i$: radial invariance holds if and only if
both $\ker\widehat\omega_i$ and $\ker\widehat\eta_i$ are
independent of $i$. By \eqref{eq:cpliftkernels}, this is exactly
\eqref{eq:cpalignment}. This proves the equivalence.
Assume henceforth that these conditions hold, so that
$TM=\E\oplus\F$ by \eqref{eq:cpsplitting}. The common
characteristic distributions of the lift are
\[
  \widehat\E=\ker\widehat\eta_i=\pi^*\E\oplus\R\partial_r,
  \qquad
  \widehat\F=\ker\widehat\omega_i=\pi^*\F\oplus\R\partial_s,
\]
where pullback subbundles are identified with their horizontal
lifts. By Theorem~\ref{thm:spalignment}, they are orthogonal,
quaternionic and parallel for $\widehat g$; in particular their ranks
are divisible by $4$. By the rank--nullity theorem,
\[
  \rank\widehat\E=\dim\widehat M-\rank\widehat\eta_i=2p_i+2,
  \qquad
  \rank\widehat\F=\dim\widehat M-\rank\widehat\omega_i=2q_i+2 .
\]
Both ranks are positive, since $\widehat\E\ni\partial_r$ and
$\widehat\F\ni\partial_s$. Hence $2p_i+2=4(l+1)$ and
$2q_i+2=4(m+1)$ for some $l,m\ge0$ independent of $i$, that is,
\[
  p_i=2l+1,\quad q_i=2m+1,
  \quad
  \rank\E=2p_i+1=4l+3,
  \quad \rank\F=2q_i+1=4m+3.
\]
Both distributions are integrable. For $\F=\charac(\alpha_1)$ this
is the computation in the proof of Theorem~\ref{thm:cspalignment}:
if $X,Y\in\Gamma(\F)$, then $\alpha_1([X,Y])=-d\alpha_1(X,Y)=0$ and
$\iota_{[X,Y]}d\alpha_1=\Lie_X(\iota_Yd\alpha_1)
-\iota_Y(\Lie_Xd\alpha_1)=0$, because $\Lie_Xd\alpha_1
=d\iota_Xd\alpha_1=0$. The same computation with $\beta_1$ applies
to $\E=\charac(\beta_1)$. Since $\E\cap\F=\{0\}$, the Frobenius
theorem gives, after shrinking, a local product
\[
  U\cong P^{4l+3}\times Q^{4m+3},
  \qquad \E|_U=TP,\qquad \F|_U=TQ,
\]
with connected factors. For $Y\in\Gamma(\F)$ and $X\in\Gamma(\E)$,
\[
  \alpha_i(Y)=0,\quad
  \Lie_Y\alpha_i=d(\alpha_i(Y))+\iota_Yd\alpha_i=0,
  \quad
  \beta_i(X)=0,\quad
  \Lie_X\beta_i=d(\beta_i(X))+\iota_Xd\beta_i=0 .
\]
As in the proof of Theorem~\ref{thm:cspalignment}, it follows that
$\alpha_i=\pr_P^*\alpha_i^P$ and $\beta_i=\pr_Q^*\beta_i^Q$ for
$1$-forms $\alpha_i^P$ on $P$ and $\beta_i^Q$ on $Q$. The volume form
\eqref{eq:cpordinary} of the $i$th pair is, up to sign, the exterior
product of $\pr_P^*\bigl(\alpha_i^P\wedge(d\alpha_i^P)^{2l+1}\bigr)$ and
$\pr_Q^*\bigl(\beta_i^Q\wedge(d\beta_i^Q)^{2m+1}\bigr)$, whose degrees
are $\dim P$ and $\dim Q$. Both factors are therefore nowhere zero,
and $\alpha_i^P$, $\beta_i^Q$ are contact forms.
Under the induced identification
\[
  U\times\R_{>0}^2
  \cong(P\times\R_{>0})\times(Q\times\R_{>0}),
\]
where $r$ is the coordinate of the first $\R_{>0}$ and $s$ that of
the second, the common distributions are
$\widehat\E=T(P\times\R_{>0})$ and $\widehat\F=T(Q\times\R_{>0})$,
and the two families of lifted forms are pulled back from these
factors, where they read
\[
  \Omega_i^P=d\!\left(\frac{r^2}{2}\alpha_i^P\right),
  \qquad
  \Omega_i^Q=d\!\left(\frac{s^2}{2}\beta_i^Q\right).
\]
Theorem~\ref{thm:spalignment} gives hyper-K\"ahler structures
$(\widehat g_P,J_i)$ on $P\times\R_{>0}$ and $(\widehat g_Q,K_i)$ on
$Q\times\R_{>0}$, with fundamental forms $\Omega_i^P$ and $\Omega_i^Q$,
and
\[
  \widehat g=\widehat g_P+\widehat g_Q,\qquad I_i=J_i\oplus K_i .
\]
Since $\Lie_R\Omega_i^P=2\Omega_i^P$ and $\Lie_S\Omega_i^Q=2\Omega_i^Q$,
the computation at the beginning of the proof of
Lemma~\ref{lem:cone} gives $\Lie_RJ_i=0$ and $\Lie_SK_i=0$. Thus both
factors are radially invariant hyper-K\"ahler symplectizations.
\end{proof}

\begin{theorem}\label{thm:cpnormalized}
A normalized hyper contact pair $(\alpha_i,\beta_i)_{i=1,2,3}$ is
locally, with its canonical metric, a Riemannian product
\begin{equation}\label{eq:cpproduct}
 (M,g)=(P^{4l+3},g_P)\times(Q^{4m+3},g_Q)
\end{equation}
of two $3$-Sasakian manifolds, and $\alpha_i$, $\beta_i$ are the
pullbacks of their respective contact forms. On the double
symplectization,
\begin{equation}\label{eq:doubleconemetric}
 \widehat g=dr\otimes dr+r^2g_P+ds\otimes ds+s^2g_Q.
\end{equation}
Conversely, every such local product, with the forms globally
defined, is a normalized hyper contact pair.
\end{theorem}

\begin{proof}
Let $\widehat g$ be the canonical metric of the lifted hyper
symplectic pair, and let $g$ be its restriction to $r=s=1$. Since
the hyper contact pair is radially invariant,
Theorem~\ref{thm:cpalignment} applies: every point of $M$ has a
neighborhood
\[
  U\cong P^{4l+3}\times Q^{4m+3}
\]
with $\alpha_i=\pr_P^*\alpha_i^P$ and $\beta_i=\pr_Q^*\beta_i^Q$,
and on $(P\times\R_{>0})\times(Q\times\R_{>0})$
\[
  \widehat g=\widehat g_P+\widehat g_Q,
\]
where $\widehat g_P$ is reconstructed from the forms
$d(r^2\alpha_i^P/2)$ and $\widehat g_Q$ from the forms
$d(s^2\beta_i^Q/2)$. The vector field $R$ is tangent to the first
factor and $S$ to the second, and each is annihilated by the metric
of the other factor. Hence the normalization
\eqref{eq:cpnormalization} restricts to
  $\widehat g_P(R,\cdot)=r\,dr$ and 
  $\widehat g_Q(S,\cdot)=s\,ds$.
Lemma~\ref{lem:cone}, applied separately to the two factors, gives
\[
  \widehat g_P=dr\otimes dr+r^2g_P,
  \qquad
  \widehat g_Q=ds\otimes ds+s^2g_Q,
\]
where $(P,g_P)$ and $(Q,g_Q)$ are $3$-Sasakian with defining
contact forms $\alpha_i^P$ and $\beta_i^Q$, respectively. This
proves \eqref{eq:doubleconemetric}, and restriction to $r=s=1$
gives $g=g_P+g_Q$, which is \eqref{eq:cpproduct}.

Conversely, suppose that the globally defined forms are locally
those of a product of two $3$-Sasakian manifolds as in the
statement. The form
\[
  \alpha_i\wedge(d\alpha_i)^{2l+1}
  \wedge\beta_i\wedge(d\beta_i)^{2m+1}
\]
is a volume form on the product, being the exterior product of the
contact volume forms of the two factors, and
$(d\alpha_i)^{2l+2}=(d\beta_i)^{2m+2}=0$. Thus each
$(\alpha_i,\beta_i)$ is a contact pair of type $(2l+1,2m+1)$.
The cones of the two $3$-Sasakian factors are hyper-K\"ahler, with
metrics $dr\otimes dr+r^2g_P$ and $ds\otimes ds+s^2g_Q$ and
fundamental forms $d(r^2\alpha_i^P/2)$ and $d(s^2\beta_i^Q/2)$,
respectively (see \cite{BG99}). The converse part of
Theorem~\ref{thm:spalignment} shows that these forms, pulled back
to the product of the two cones, define a hyper symplectic pair
with canonical metric \eqref{eq:doubleconemetric} and fundamental
forms \eqref{eq:cpforms}. By Definition~\ref{def:cp}, the original
pairs form a hyper contact pair. Their characteristic distributions
$\charac(\alpha_i)=TQ$ and $\charac(\beta_i)=TP$ do not depend on
$i$, so Theorem~\ref{thm:cpalignment} gives radial invariance.
For any tangent vector $V$ on the double symplectization,
\[
  \widehat g(R,V)=dr(R)\,dr(V)=r\,dr(V),
  \qquad
  \widehat g(S,V)=ds(S)\,ds(V)=s\,ds(V),
\]
because the pullbacks of $g_P$ and $g_Q$ annihilate the radial
vectors and $ds(R)=dr(S)=0$. Hence the pair is normalized. Finally,
the metric and the endomorphisms reconstructed from the forms are
uniquely determined by them, so the local constructions agree on
overlaps, and the conditions above hold globally.
\end{proof}

\begin{proposition}\label{prop:Ftensors}
Let $(\alpha_i,\beta_i)_{i=1,2,3}$ be a normalized hyper contact
pair on $M$, with canonical metric $g$. Let $Z_i$, $W_i$ be the
vector fields and $F_i$ the endomorphisms of $TM$ defined by
\[
 g(Z_i,\cdot)=\alpha_i,\qquad g(W_i,\cdot)=\beta_i,\qquad
 g(\cdot,F_i\,\cdot)=\frac12\,(d\alpha_i+d\beta_i).
\]
Then $Z_i,W_i$ are the Reeb fields of $(\alpha_i,\beta_i)$, the six
fields $Z_1,Z_2,Z_3,W_1,W_2,W_3$ are orthonormal, and for every
cyclic $(i,j,k)$
\begin{gather}
 F_i^2=-\Id+\alpha_i\otimes Z_i+\beta_i\otimes W_i,\label{eq:Fsquare}\\
 F_iF_j=F_k+\alpha_j\otimes Z_i+\beta_j\otimes W_i,\label{eq:Fquaternion}\\
 F_iZ_j=Z_k,\quad F_iW_j=W_k,\label{eq:FReeb}\\
 d\alpha_i+d\beta_i=2g(\cdot,F_i\cdot).\label{eq:associatedmetric}
\end{gather}
Moreover, for every cyclic $(i,j,k)$,
\begin{equation}\label{eq:Reebbrackets}
 [Z_i,Z_j]=2Z_k,\qquad [W_i,W_j]=2W_k,
\end{equation}
and, for all $a,b\in\{1,2,3\}$, $[Z_a,W_b]=0$. Furthermore, for every
cyclic $(i,j,k)$,
\begin{equation}\label{eq:alphaown}
\alpha_i=\frac12\iota_{Z_j}d\alpha_k
          =-\frac12\iota_{Z_k}d\alpha_j, \qquad
 \beta_i=\frac12\iota_{W_j}d\beta_k
          =-\frac12\iota_{W_k}d\beta_j,
\end{equation}
and, for all $a,b\in\{1,2,3\}$,
\begin{equation}\label{eq:crosszero}
 \iota_{Z_a}d\beta_b=\iota_{W_a}d\alpha_b=0.
\end{equation}
Finally, if $\pi_{\E},\pi_{\F}$ are the projections associated with
$TM=\E\oplus\F$, then
\begin{align}\label{eq:nablaF}
 (\nabla_XF_i)Y
 ={}&g(\pi_{\E}X,\pi_{\E}Y)Z_i
     -\alpha_i(Y)\pi_{\E}X+g(\pi_{\F}X,\pi_{\F}Y)W_i
     -\beta_i(Y)\pi_{\F}X.
\end{align}
\end{proposition}

\begin{proof}
The tensors $Z_i,W_i,F_i$ are globally defined, since $g$ is
nondegenerate. Let $(\phi_i^P,\xi_i^P,\alpha_i^P,g_P)$ and
$(\phi_i^Q,\xi_i^Q,\beta_i^Q,g_Q)$ be the $3$-Sasakian structures of
the local factors in Theorem~\ref{thm:cpnormalized}. Since
$g=g_P+g_Q$, $d\alpha_i^P=2g_P(\cdot,\phi_i^P\cdot)$ and
$d\beta_i^Q=2g_Q(\cdot,\phi_i^Q\cdot)$, the defining identities give
the local expressions
\[
 Z_i=\xi_i^P,\qquad W_i=\xi_i^Q,\qquad F_i=\phi_i^P\oplus\phi_i^Q .
\]
Equivalently, on $r=s=1$ these tensors are recovered from the
hyper-K\"ahler structure of the double cone by
\[
 Z_i=-I_iR,\qquad W_i=-I_iS,\qquad
 F_iX=I_iX-\alpha_i(X)R-\beta_i(X)S .
\]
Since $\alpha_i$ comes from $P$ and $\beta_i$ from $Q$, the
fields $\xi_i^P,\xi_i^Q$ satisfy the defining conditions of the
Reeb fields of $(\alpha_i,\beta_i)$, so $Z_i,W_i$ are those Reeb
fields. They are orthonormal because the $\xi_i^P$ are orthonormal
in $(P,g_P)$, the $\xi_i^Q$ are orthonormal in $(Q,g_Q)$, and the
two factors are orthogonal.
The common projector of the lifted hyper symplectic pair, given by
Theorem~\ref{thm:spalignment}, restricts on $TM$ at $r=s=1$ to
$\pi_{\E}$, and the complementary projector restricts to $\pi_{\F}$.

The almost contact metric $3$-structure identities on each factor
give \eqref{eq:Fsquare}--\eqref{eq:FReeb}; for example, the tangential
component of $I_iI_jX=I_kX$ is \eqref{eq:Fquaternion}. Identity
\eqref{eq:associatedmetric} holds by definition. The Sasakian
identity $\nabla^P_X\xi_i^P=-\phi_i^PX$ gives, as in the proof of
Theorem~\ref{thm:cspnormalized},
$[\xi_i^P,\xi_j^P]=-\phi_j^P\xi_i^P+\phi_i^P\xi_j^P=2\xi_k^P$, and
similarly on $Q$; this proves \eqref{eq:Reebbrackets}. Vector fields
lifted from different factors commute, so $[Z_a,W_b]=0$. The
computation proving \eqref{eq:correctcontractions} in
Theorem~\ref{thm:cspnormalized}, carried out on each factor, gives
\eqref{eq:alphaown}. Identity
\eqref{eq:crosszero} holds because $d\beta_b$ is pulled back from
$Q$ while $Z_a$ is tangent to $P$, and $d\alpha_b$ is pulled back
from $P$ while $W_a$ is tangent to $Q$. Finally, the Levi-Civita
connection of $g=g_P+g_Q$ is the product connection; it preserves
$\E$ and $\F$, and applying \eqref{eq:sasakiconvention} on each
factor gives \eqref{eq:nablaF}.
\end{proof}

We now use foliated bundles to construct compact examples for which the canonical local product decomposition does not extend to a global product. As we will see, a special phenomenon appears in the case of hyper contact pairs. Keeping the notation of Section~\ref{sec:csp}, we have the following Proposition whose proof is similar to that of Proposition~\ref{prop:cspfoliated}.

\begin{proposition}\label{prop:cpfoliated}
Let $(B^{4l+3},g_B)$ and $(F^{4m+3},g_F)$ be $3$-Sasakian manifolds,
with contact forms $\alpha_i^B$ and $\beta_i^F$, and let
$\rho:\pi_1(B)\to\operatorname{Aut}(F)$ be a homomorphism. Then the
forms $\pr_1^*\alpha_i^{\widetilde B}$ and $\pr_2^*\beta_i^F$ on
$\widetilde B\times F$, where $\alpha_i^{\widetilde B}$ is the pullback
of $\alpha_i^B$, descend to $1$-forms $\alpha_i,\beta_i$ on
$M=\widetilde B\times_\rho F$, and $(\alpha_i,\beta_i)_{i=1,2,3}$ is a
normalized hyper contact pair. Its canonical metric is induced by
$g_{\widetilde B}+g_F$. The distribution $\charac(\alpha_i)$ is
tangent to the fibres, and $\charac(\beta_i)$ is tangent to the
foliation defined by the slices $\widetilde B\times\{y\}$. If $B$ and
$F$ are compact, so is $M$.
\end{proposition}

\begin{example}\label{ex:cpspheres}
Let $\Gamma\subset\Sp(1)$ be a nontrivial finite subgroup, for
instance a cyclic group $\Z_k$ or a binary polyhedral group. Let
$\Gamma$ act on $\HH^{n+1}$ by left multiplication by $q\cdot\Id$,
$q\in\Gamma$. This action lies in $\Sp(n+1)$, and it is free on
$S^{4n+3}$: if $qx=x$ with $x\neq0$, then $q=1$. Take
\[
 B=S^{4l+3}/\Gamma,\qquad F=S^{4m+3},\qquad
 \rho:\Gamma\hookrightarrow\Sp(1)\subset\Sp(m+1),
\]
so that $\widetilde B=S^{4l+3}$ and $\pi_1(B)=\Gamma$. By
Proposition~\ref{prop:cpfoliated},
\[
 M=(S^{4l+3}\times S^{4m+3})/\Gamma,
 \qquad q\cdot(x,y)=(qx,qy),
\]
is a compact manifold of dimension $4(l+m)+6$ with a normalized
hyper contact pair. The smallest case is $l=m=0$, which gives
the six-dimensional manifolds $(S^3\times S^3)/\Gamma$.

The splitting $TM=\E\oplus\F$ is not a global Riemannian product.
Indeed, the leaf of $\E$ through $[x,y]$ is the image of
$S^{4l+3}\times\{y\}$, and the leaf of $\F$ through $[x',y']$ is the
image of $\{x'\}\times S^{4m+3}$. They meet exactly at the points
$[qx',y]$, $q\in\Gamma$, and these are pairwise distinct, because
$\Gamma$ acts freely on both spheres. Hence every leaf of $\E$ meets
every leaf of $\F$ in exactly $|\Gamma|>1$ points, whereas in a
Riemannian product $P\times Q$ the leaves $P\times\{y\}$ and
$\{x\}\times Q$ meet in one point. 
\end{example}

Every compact connected manifold carrying a normalized hyper contact pair is
a finite diagonal quotient of a product of compact simply connected 3-Sasakian
manifolds.

\begin{proposition}\label{prop:cpcompact}
Let $(\alpha_i,\beta_i)$ be a normalized hyper contact pair on a
compact connected manifold $M$, with canonical metric $g$. Then
$\pi_1(M)$ is finite, and
\[
 M\cong(\widetilde P\times\widetilde Q)/\Gamma,
 \qquad \Gamma\cong\pi_1(M),
\]
where $\widetilde P$ and $\widetilde Q$ are compact simply connected
$3$-Sasakian manifolds and $\Gamma$ is a finite subgroup of
$\operatorname{Aut}(\widetilde P)\times\operatorname{Aut}(\widetilde Q)$
acting freely on $\widetilde P\times\widetilde Q$. The
diffeomorphism identifies $g$ with the metric induced by
$g_{\widetilde P}+g_{\widetilde Q}$, and $\alpha_i$, $\beta_i$ with
the forms induced by the contact forms of $\widetilde P$ and of
$\widetilde Q$, respectively. If the projection
$\Gamma\to\operatorname{Aut}(\widetilde P)$ is injective and its
image $\Gamma_P$ acts freely on $\widetilde P$, then $M$ is the
foliated bundle $\widetilde P\times_\rho\widetilde Q$ over
$\widetilde P/\Gamma_P$, where $\rho$ is the second projection.
\end{proposition}

\begin{proof}
We proceed in several steps.

By Theorem~\ref{thm:cpalignment}, $TM=\E\oplus\F$ with
$\E=\charac(\beta_i)$ and $\F=\charac(\alpha_i)$ for every $i$. By
Theorem~\ref{thm:cpnormalized}, every point of $M$ has a
neighborhood $U\cong P\times Q$ on which $g=g_P+g_Q$, and 
$\E=TP$ and $\F=TQ$. Thus $\E$ and $\F$ are $g$-orthogonal. The
Levi-Civita connection of a Riemannian product is the product
connection, which preserves the tangent bundles of the two factors.

Let $\pi:\widetilde M\to M$ be the universal covering, and let
$\widetilde g$, $\widetilde\alpha_i$, $\widetilde\beta_i$ be the
pullbacks of $g$, $\alpha_i$, $\beta_i$. All the local statements of
Theorems~\ref{thm:cpalignment} and~\ref{thm:cpnormalized} hold on
$\widetilde M$, since $\pi$ is a local isometry preserving the forms.
In particular, $T\widetilde M=\widetilde\E\oplus\widetilde\F$ with
$\widetilde\E=\pi^*\E$ and $\widetilde\F=\pi^*\F$ orthogonal and
parallel. Since $M$ is compact, $g$ is complete, and so is
$\widetilde g$. By the de Rham decomposition theorem, applied to the
complete simply connected manifold $(\widetilde M,\widetilde g)$ and
to the parallel orthogonal splitting
$\widetilde\E\oplus\widetilde\F$, there is an isometry
\[
 \Phi:(\widetilde P\times\widetilde Q,\,
 g_{\widetilde P}+g_{\widetilde Q})\longrightarrow
 (\widetilde M,\widetilde g),
\]
where $\widetilde P$ and $\widetilde Q$ are leaves of
$\widetilde\E$ and of $\widetilde\F$, with their induced metrics,
and $d\Phi$ maps $T\widetilde P\oplus0$ onto $\widetilde\E$ and
$0\oplus T\widetilde Q$ onto $\widetilde\F$. We identify
$\widetilde M$ with $\widetilde P\times\widetilde Q$ via $\Phi$. The
leaves of $\widetilde\E$ are the slices $\widetilde P\times\{q\}$,
the leaves of $\widetilde\F$ are the slices $\{p\}\times\widetilde Q$,
and both factors are complete. Since
$\pi_1(\widetilde M)=1$, both factors are simply connected.

Let $X$ be a section of $\widetilde\F=\charac(\widetilde\alpha_i)$. Then
$\iota_X\widetilde\alpha_i=0$ and $\iota_Xd\widetilde\alpha_i=0$, so
$\Lie_X\widetilde\alpha_i=0$ by Cartan's formula. Hence
$\widetilde\alpha_i$ is basic for the submersion
$\pr_{\widetilde P}:\widetilde P\times\widetilde Q\to\widetilde P$,
whose fibres $\{p\}\times\widetilde Q$ are connected. Thus
$\widetilde\alpha_i=\pr_{\widetilde P}^*a_i$, where
$a_i=\iota_q^*\widetilde\alpha_i$ for the inclusion
$\iota_q:\widetilde P\to\widetilde P\times\{q\}$, and $a_i$ does not
depend on $q$. In the same way,
$\widetilde\beta_i=\pr_{\widetilde Q}^*b_i$, with $b_i$ a $1$-form on
$\widetilde Q$. 
Therefore $(\widetilde P,g_{\widetilde P})$,
with the forms $a_i$, is $3$-Sasakian of dimension $4l+3$ with contact forms $a_i$.
Likewise, $\widetilde Q$ is $3$-Sasakian of dimension $4m+3$ with
contact forms $b_i$.

A $3$-Sasakian manifold of dimension $4n+3$ is Einstein with
$\Ric=(4n+2)g$ (see \cite{BG99}). Thus $\widetilde P$ is complete with
$\Ric_{\widetilde P}=(4l+2)g_{\widetilde P}
=(\dim\widetilde P-1)g_{\widetilde P}$, and $4l+2>0$. By Myers'
theorem $\widetilde P$ is compact, with diameter at most $\pi$. The
same argument shows that $\widetilde Q$ is compact, and so is
$\widetilde M=\widetilde P\times\widetilde Q$. Therefore $\pi_1(M)$ is finite.

Let $\Gamma\cong\pi_1(M)$ be the group of deck transformations, and let
$\gamma\in\Gamma$. Then $\gamma$ is an isometry of $\widetilde g$
with $\gamma^*\widetilde\alpha_i=\widetilde\alpha_i$ and
$\gamma^*\widetilde\beta_i=\widetilde\beta_i$, because these tensors
are pulled back from $M$. Hence $d\gamma$ preserves
$\widetilde\E=\charac(\widetilde\beta_i)$ and
$\widetilde\F=\charac(\widetilde\alpha_i)$, and $\gamma$ maps leaves
to leaves.
Consequently, $\pr_{\widetilde Q}\circ\gamma(p,q)$ depends only on
$q$, and $\pr_{\widetilde P}\circ\gamma(p,q)$ depends only on $p$.
Fixing $(p_0,q_0)$ and setting
\[
 f_P(p)=\pr_{\widetilde P}\gamma(p,q_0),\qquad
 f_Q(q)=\pr_{\widetilde Q}\gamma(p_0,q),
\]
we obtain smooth maps with $\gamma=f_P\times f_Q$. Applying the same
argument to $\gamma^{-1}$ shows that $f_P$ and $f_Q$ are
diffeomorphisms. Since $d\gamma=df_P\oplus df_Q$ is orthogonal for
$g_{\widetilde P}+g_{\widetilde Q}$, both $f_P$ and $f_Q$ are
isometries. Moreover,
\[
 \pr_{\widetilde P}^*a_i=\widetilde\alpha_i
 =\gamma^*\widetilde\alpha_i
 =\gamma^*\pr_{\widetilde P}^*a_i
 =\pr_{\widetilde P}^*f_P^*a_i ,
\]
and $\pr_{\widetilde P}$ is a surjective submersion, so
$f_P^*a_i=a_i$. Similarly $f_Q^*b_i=b_i$. Thus
$f_P\in\operatorname{Aut}(\widetilde P)$ and
$f_Q\in\operatorname{Aut}(\widetilde Q)$. The map
$\gamma\mapsto(f_P,f_Q)$ is an injective group homomorphism
$\Gamma\to\operatorname{Aut}(\widetilde P)\times
\operatorname{Aut}(\widetilde Q)$, and its image acts freely on
$\widetilde P\times\widetilde Q$. Finally,
$M=\widetilde M/\Gamma=(\widetilde P\times\widetilde Q)/\Gamma$, and
$g$, $\alpha_i$, $\beta_i$ are the tensors induced by
$g_{\widetilde P}+g_{\widetilde Q}$, $\pr_{\widetilde P}^*a_i$ and
$\pr_{\widetilde Q}^*b_i$.

Now assume that $h:\Gamma\to\operatorname{Aut}(\widetilde P)$,
$\gamma\mapsto f_P$, is injective and that its image $\Gamma_P$ acts
freely on $\widetilde P$. The finite group $\Gamma_P$ acts freely by
isometries preserving the $a_i$. Hence $B=\widetilde P/\Gamma_P$ is a
compact $3$-Sasakian manifold. Its universal covering is
$\widetilde P\to B$, because $\widetilde P$ is simply connected, and
$\pi_1(B)\cong\Gamma_P$. Let
$\rho=\pr_2\circ h^{-1}:\Gamma_P\to\operatorname{Aut}(\widetilde Q)$,
so that $h^{-1}(f)=f\times\rho(f)$ for every $f\in\Gamma_P$. Then the
action of $\Gamma$ on $\widetilde P\times\widetilde Q$ is
\[
 f\cdot(p,q)=(f(p),\rho(f)(q)),\qquad f\in\Gamma_P=\pi_1(B),
\]
which is the action defining $\widetilde P\times_\rho\widetilde Q$.
Thus $M\cong\widetilde P\times_\rho\widetilde Q$, and by the previous arguments the
hyper contact pair of $M$ is the one given by
Proposition~\ref{prop:cpfoliated}.
\end{proof}

\begin{remark}\label{rem:cpvscsp}
The finite-fundamental-group conclusion of Proposition~\ref{prop:cpcompact} fails for normalized
hyper contact–symplectic pairs. In Example~\ref{ex:cspinfinite} the fundamental group is $\Z^{4m}$, the suspension holonomy is infinite, and the leaves of $\F$ are noncompact. For normalized hyper contact pairs both local factors are $3$-Sasakian and have positive Ricci curvature. For normalized hyper contact–symplectic pairs the second factor
is hyper-Kähler and Ricci-flat.
\end{remark}

\bigskip
\noindent\textbf{AI Disclosure.} Claude Opus 5.5 and ChatGPT-6 Astra were
used in the preparation of this manuscript. All mathematical ideas and results are
due to the authors, who take full responsibility for the mathematical
content and the final text.

\end{document}